\documentclass{article}

\PassOptionsToPackage{numbers, compress}{natbib}
\usepackage[final]{neurips_2020}

\usepackage{times}
\usepackage[utf8]{inputenc} 
\usepackage[T1]{fontenc}    
\usepackage{hyperref}       
\usepackage{url}            
\usepackage{booktabs}       
\usepackage{amsfonts}       
\usepackage{nicefrac}       
\usepackage{microtype}      

\usepackage{latex-macros}
\usepackage[labelfont=bf]{caption}
\usepackage{xfrac}
\usepackage{multirow}
\newcommand{\gap}{\mathsf{gap}}
\newcommand{\sep}{\mathsf{sep}}

\newcommand{\twoinf}{\ell_{2 \to \infty}}

\newtheorem{theorem}{Theorem}
\newtheorem{lemma}{Lemma}

\newtheorem{proposition}{Proposition}
\newtheorem{corollary}{Corollary}
\newtheorem{definition}{Definition}
\newtheorem{assumption}{Assumption}
\newtheorem*{myprop}{Proposition}

\let\counterwithin\relax
\usepackage{chngcntr}
\usepackage{apptools}
\AtAppendix{\counterwithin{lemma}{section}}
\AtAppendix{\counterwithin{theorem}{section}}

\theoremstyle{remark}

\usepackage{xcolor}
\usepackage[capitalize,noabbrev]{cleveref}
\usepackage{tikz,pgfplots}
\pgfplotsset{compat=1.14}
\usepackage{cancel}
\usepackage{algorithm,algpseudocode}
\definecolor{lblue}{HTML}{908cc0}
\definecolor{mblue}{HTML}{519cc8}
\definecolor{hblue}{HTML}{1d5996}
\definecolor{lred}{HTML}{cb5501}
\definecolor{mred}{HTML}{f1885b}
\definecolor{hred}{HTML}{b3001e}
\definecolor{ttred}{HTML}{ca3542}

\usepackage{svg}

\usepackage{hyperref}
\definecolor{mylinkcolor}{RGB}{0,0,140}
\hypersetup{colorlinks,allcolors=mylinkcolor,citecolor=mylinkcolor}

\title{Entrywise Convergence of \\ Iterative Methods for Eigenproblems}

\author{%
  Vasileios Charisopoulos\\
  Department of Operations Research \& Information Engineering\\
  Cornell University\\
  Ithaca, NY 14853\\
  \texttt{vc333@cornell.edu}
  \And
  Austin R.~Benson\\
  Department of Computer Science\\
  Cornell University\\
  Ithaca, NY 14853\\
  \texttt{arb@cs.cornell.edu}
  \And
  Anil Damle\\
  Department of Computer Science\\
  Cornell University\\
  Ithaca, NY 14853\\
  \texttt{damle@cornell.edu}
}

\begin{document}

\maketitle

\begin{abstract}
Several problems in machine learning, statistics, and other fields rely on
computing eigenvectors. For large scale problems, the computation of these eigenvectors is typically performed via iterative schemes
such as subspace iteration or Krylov methods. While there is classical and comprehensive analysis for subspace convergence guarantees with respect to the spectral norm, in many modern applications other notions
of subspace distance are more appropriate. Recent theoretical
work has focused on perturbations of subspaces measured in the $\twoinf$ norm, but does not consider the actual computation of eigenvectors. Here we address the convergence of subspace iteration
when distances are measured in the $\twoinf$ norm and provide deterministic bounds. We complement our analysis with a practical stopping criterion and demonstrate its applicability via numerical experiments. Our results show that one can get comparable performance on downstream tasks while requiring fewer iterations, thereby saving substantial computational time.
\end{abstract}

\section{Introduction \& Background}
\label{sec:intro}
Spectral methods play a fundamental role in machine learning, statistics, and data mining.
Methods for foundational tasks such as
clustering~\cite{VonLux07};
semi-supervised learning~\cite{mahoney2012local};
dimensionality reduction~\cite{belkin2002laplacian,friedman2001elements,roweis2000nonlinear};
latent factor models~\cite{gower2014netflix}
ranking and preference learning~\cite{maystre2015fast,vigna2016spectral}; graph
signal processing~\cite{OFK+18,SNF+13}; and covariance estimation
all use information about eigenvalues and eigenvectors (or singular values and singular vectors)
from an underlying data matrix (either directly or indirectly).
The pervasiveness of spectral methods in machine learning applications\footnote{
    For example, searching for ``\texttt{arpack}''~\cite{arpack} (an iterative eigensolver) in
    the \texttt{scikit-learn}~\cite{sklearn} Github repository reveals that several modules
    depend on it crucially.}
has greatly influenced the last decade of research in large-scale
computation, including but not limited to sketching / randomized NLA~\cite{HMT11,MT20,Woodruff14}
as well as theoretical guarantees for linear algebra primitives (\emph{e.g.}, eigensolvers,
low-rank decompositions) in previously overlooked settings.
%

In many of these cases, the relevant information is in the ``leading''
eigenvectors, \emph{i.e.,} those corresponding to the $k$ algebraically largest eigenvalues
for some $k$ (possibly after shifting and rescaling). To avoid performing a full
eigendecomposition, these are typically approximated with iterative algorithms such
as the power or Lanczos methods.
The approximation quality, as measured by subspace distance (equivalent to using the $\ell_2$ norm,
up to rotation), is well-understood and enjoys comprehensive convergence analysis~\cite{Demmel97,GVL13,Parlett98,Saad11}.

While spectral norm error analysis has been the standard-bearer for numerical
analysis, recent work has considered different subspace distance
measures~\cite{CLC+19,CCF18,FWZZ18,XiaZhou19}.
The motivation for these changes is statistical, as opposed to numerical:
we observe a matrix $\tilde{A} = A + E$,
where $E$ is a source of noise and $A = \expec{\tilde{A}}$ is the ``population''
version of $A$, containing the desired spectral information. We are then
interested in $\norm{\tilde{u}_i \pm u_i}_{\infty}$ as a distance measure between the
eigenvectors of $\tilde{A}$ and $A$.  Here, the $\ell_{\infty}$ norm captures
``entry-wise'' error and is more appropriate when we care about maximum
deviation; for example, when entries of the eigenvector are used to rank nodes
or provide cluster assignments in a graph. This type of distance is often much
smaller than the spectral norm and, in contrast to the latter, reveals information
about the distribution of error over the entries.
Recent theoretical results relate the noise $E$ to the perturbation
in the eigenvectors, as measured by $\ell_{\infty}$ or $\twoinf$ norm errors~\cite{CapeTangPriebe19,DamleSun19,FanWangZhong18,KolXia16,LL19}.
Moreover, these results are often directly connected to machine learning
problems~\cite{AFWZ17,CFMW19,EldBelWang18,ZB18}.

The message from this body of literature is that when eigenvectors are interpreted entry-wise, we should
measure our error entry-wise as well. The aforementioned works show what we can do
\emph{if} we have eigenvectors satisfying perturbation bounds in a different norm,
but do not address their \emph{computation}. Numerical algorithms typically use
the $\ell_2$ norm, yet the motivation for norms like $\twoinf$ is that $\ell_2$ can
be a severe overestimate for the relevant approximation quality.
Moreover, despite the long history of research into stopping criteria
for iterative methods in the unitarily-invariant
setting~\cite{ADR92,BDM93,BennaniBraconnier94,GolubMeurant97,LehSorYang97},
there are no generic stopping criteria closely tracking the
quality of an approximation in the $\twoinf$ norm. For example,
downstream tasks that depend on entrywise ordering, such as graph bipartitioning
via the (approximate) Fiedler vector~\cite{fairbanks2016new} or
spectral ranking via the \textit{Katz} centrality~\cite{NSF+17} employ $\ell_2$
bounds, when instead the $\ell_{\infty}$ norm would constitute a better proxy.
Some local spectral graph partitioning methods can be written as iteratively
approximating an eigenvector in a (scaled) $\ell_{\infty}$
norm~\cite{andersen2006local}, but these algorithms are far more specialized
than general eigensolvers.
The situation is similar when using more than one eigenvector; in spectral
clustering with $r$ clusters, after an appropriate rotation of the eigenvector
matrix, the magnitude of the elements in the $i^{\text{th}}$ row measures the per-cluster
membership likelihood of the $i^{\text{th}}$ node, making the $\ell_{2 \to \infty}$ norm (which is invariant to
unitary transformations on the right) a more appropriate distance measure than the
spectral norm (see \textit{e.g.}, \cite{LST+14}).

Here, we bridge this gap by providing an analysis for the convergence of subspace iteration, a
widely-used iterative method for computing leading eigenvectors, in terms of
$\twoinf$ errors. We complement that with a practical stopping criterion applicable to any iterative method for invariant subspace computation that
tracks the $\twoinf$ error of the approximation. Our results show how,
for a given error tolerance, one can perform many fewer subspace iterations to
get the same desired performance on a downstream task that uses the
eigenvectors (or, more generally, an invariant subspace) --- as $\norm{V}_{2\to\infty}
\in [1, \sqrt{r}] \max_{i,j}\abs{V_{ij}}$ for $V \in \Rbb^{n \times r}$, and often
$r \ll n$, our bounds are also a good ``proxy'' for the maximum entrywise error.
The aforementioned reduction in iterations directly translates to substantial reductions in
computation time.
We demonstrate our methods with the help of applications involving real-world
graph data, including node ranking in graphs, sweep cut profiles for spectral
bipartioning, and general spectral clustering.

\subsection{Notation}
We use the standard inner product on Euclidean spaces, defined by
$\ip{X, Y} := \trace{X^\T Y}$ for vectors/matrices $X, Y$. We write $\Obb_{n,
k}$ for the set of matrices $U \in \Rbb^{n \times k}$ such that $U^\T U =
I_k$, dropping the second subscript when $n = k$.
We use standard notation for norms, namely $\norm{A}_2 := \sup_{x: \norm{x}_2 = 1}
\norm{Ax}_2$ and $\norm{A}_F :=
\sqrt{\ip{A, A}}$. Moreover, we remind the reader that the $\ell_{\infty} \to
\ell_{\infty}$ operator norm for a matrix $A \in \Rbb^{m \times n}$ is given by
\(
    \norm{A}_{\infty} := \max_{i \in [m]} \norm{A_{i, :}}_1,
\)
where $A_{i, :}$ denotes the $i^{\text{th}}$ row of $A$ and $A_{:, i}$ denotes
its $i^{\text{th}}$ column. Finally, the $\ell_{2 \to \infty}$ norm is defined
by
\begin{equation}
    \norm{A}_{2 \to \infty} := \sup_{x: \norm{x}_2 = 1} \norm{Ax}_{\infty}
    =
 	\max_{i \in [m]} \norm{A_{i, :}}_2.
    \label{eq:2-infty-norm}
\end{equation}

\textbf{Subspace distances.} Given two orthogonal matrices $V, \tilde{V} \in \Obb_{n, r}$ inducing subspaces
$\cV, \tilde{\cV}$, their so-called subspace distance is defined as
\(
    \dist_{2}(V, \tilde{V}) := \|VV^\T-\tilde{V}\tilde{V}^\T\|_2,
\)
with several equivalent definitions, \emph{e.g.}, via the concept of \textit{principal
angles}, or via $\norm{V_{\perp}^\T \tilde{V}}_2$, where $V_{\perp}$ is a basis for the
subspace orthogonal to $\cV$. Here we will use a slightly different notion of distance between subspaces with respect to $\norm{\cdot}_{2\to \infty}$ defined as
\begin{equation}
    \dist_{2\to\infty}(V, \tilde{V}) := \inf_{O \in \Obb_{r, r}}
    \norm{V - \tilde{V} O}_{2\to \infty}.
    \label{eq:dist-2-infty}
\end{equation}
This metric allows us to control errors in a ``row-wise'' or ``entry-wise''
sense; for example, in the case where $r=1$ this reduces to infinity norm
control over the differences between eigenvectors.
Finally, some of the stated results use the \textit{separation between matrices}
measured along a linear subspace (with respect to some norm $\norm{\cdot}_{\star}$):
\begin{equation}
    \sep_{\star, W}(B, C) = \inf\set{\norm{ZB - CZ}_{\star} \mmid
        \norm{Z}_{\star} = 1, Z \in \mathrm{range}(W)}
    \label{eq:sep}
\end{equation}
When $\norm{\cdot}_{\star}$ is unitarily invariant and $B, C$ are diagonal, we recover
$ \sep_{\star, W}(B, WCW^\T) = \lambda_{\min}(B) - \lambda_{\max}(C) $;
thus $\sep$ generalizes the notion of an eigengap.

\section{Convergence of subspace iteration}
\label{sec:convergence}
In this section, we analyze the convergence of subspace iteration
(Algorithm~\ref{alg:subspace-iteration}) with respect
to the $\ell_{2 \to \infty}$ distance. In particular, we assume that we are
working with a symmetric matrix $A$ with eigenvalue decomposition
\begin{equation}
    A = V \Lambda V^\T + V_{\perp} \Lambda_{\perp} V_{\perp}^\T,
    \label{eq:eigendecomp}
\end{equation}
where $\Lambda, \Lambda_{\perp}$ are diagonal matrices containing the $r$
largest and $n - r$ smallest eigenvalues of $A$. 
For simplicity, we assume that the eigenvalues satisfy
$\lambda_1(A) \geq \dots \geq \lambda_r(A) > \lambda_{r+1}(A) \geq \dots
\lambda_n(A)$ and, furthermore, that $\min_{k =
1,\ldots,r}\lvert\lambda_k(A)\rvert>\max_{k =
r+1,\ldots,n}\lvert\lambda_k(A)\rvert$.\footnote{Our results hold for the
largest magnitude eigenvalues assuming one defines the eigenvalue gap
appropriately later. The simplification to the $r$ algebraically largest
eigenvalues being the largest in magnitude avoids burdensome notation without
losing anything essential.}

Our perturbation bounds and stopping criterion both involve the \textit{coherence}
of the principal eigenvector matrix, which is a standard assumption in compressed
sensing~\cite{CandesRecht09}.
\begin{definition}[Coherence] \label{defn:incoherence}
    Given $V \in \Obb_{n, r}$, we define its \textbf{coherence} as the
    smallest $\mu > 0$ such that
    \begin{equation}
        \norm{V}_{2\to\infty} = \max_{i \in [n]} \norm{VV^{\T} e_i}_2 
        \leq \mu \sqrt{\frac{r}{n}}.
    \end{equation}
\end{definition}
Given~\Cref{defn:incoherence}, a matrix of eigenvectors is \textit{incoherent}
if none of its rows have a large element (i.e. all elements are on the order of
$\sqrt{1 /n}$).


\begin{algorithm}[tb]
	\caption{Subspace iteration}
	\label{alg:subspace-iteration}
	\begin{algorithmic}
		\State \textbf{Input}: initial guess $Q_0 \in \Obb_{n, k}$, symmetric matrix $A$,
		iterations $T$
		\For{$t = 1, 2, \dots, T$}
			\State $V^{(t)} := A Q_{t-1}$; \quad
            $Q_t, R_t = \texttt{qr}(V^{(t)})$ \Comment{QR decomposition}
		\EndFor
		\State \Return $Q_T$
	\end{algorithmic}
\end{algorithm}

The following result shows that $\dist_{2 \to \infty}(Q_t, V)$ can be considerably smaller than $\dist_{2}(Q_t, V)$. 
Unfortunately, our analysis involves the unwieldy term $\norm{V_{\perp} \Lambda_{\perp}^t
V_{\perp}^\T}_{\infty}$, which is nontrivial to upper bound to
obtain a better rate than that obtained using norm equivalence. To circumvent
this, we impose a technical assumption.
\begin{assumption} \label{asm:norm-bound}
    For the matrix of interest, $V_{\perp}$ satisfies
    \begin{equation}
        \norm{V_{\perp} \Lambda_{\perp}^t V_{\perp}^\T}_{\infty} \leq C \cdot
        \lambda_{r+1}^t \norm{V_{\perp} V_{\perp}^\T}_{\infty},
        \label{eq:assumption-2}
    \end{equation}
    for a small constant $C$ and all powers $t \in \Nbb$.
\end{assumption}
Assumption~\ref{asm:norm-bound} arises due to our proof technique, and may
be removed by a more careful analysis
(the supplement contains a preliminary result in this direction).
We \textit{empirically verified that it holds} with a constant $C < 2$,
for all powers $t$ up to the last elapsed iteration of Algorithm~\ref{alg:subspace-iteration} in our numerical experiments of
Section~\ref{sec:applications}; this makes us rather
confident that it is a reasonable assumption in real-world datasets.
\begin{proposition}  \label{prop:2inf-convergence-general}
    Suppose Assumption~\ref{asm:norm-bound} holds.
    The iterates $\set{Q_t}$ produced by~\Cref{alg:subspace-iteration} with
    initial guess $Q_0$ satisfy
    \begin{align}
        \begin{aligned}
        & \dist_{2 \to \infty}(Q_t, V) \leq \left(\frac{\lambda_{r +
        1}}{\lambda_r}
        \right)^t \left[ \mu\sqrt{\frac{2r}{n}} \frac{d_0}{\sqrt{1 - d_0^2}} +
        \frac{C(1 + \mu \sqrt{r})}{\sqrt{1 - d_0^2}}
        \dist_{2 \to \infty}(Q_0, V) \right],
        \end{aligned}
        \label{eq:infty-dist}
    \end{align}
    where $d_0 := \norm{Q_0^\T V_{\perp}}_2 \equiv \dist_2(Q_0, V)$,
    $r = \dim(V)$, and $\mu$ is the coherence of $V$.
\end{proposition}

When $\lambda_{r+2} \ll \lambda_{r+1}$, a slight modification of the above
proof yields a refined upper bound.
\begin{proposition}  \label{prop:2infnorm-convergence-single}
    The iterates $\set{Q_t}$ produced by~\Cref{alg:subspace-iteration} with initial guess $Q_0$ satisfy
    \begin{align}
    \label{eq:infty-dist-single}
    \begin{aligned}
        \dist_{2 \to \infty}(Q_t, V) & \leq
        \left(\frac{\lambda_{r+1}}{\lambda_r}\right)^t
        \left[ \mu \sqrt{\frac{2r}{n}}
        \cdot \frac{d_0}{\sqrt{1 - d_0^2}} +
        \frac{\norm{v_{r+1} v_{r+1}^\T}_{\infty}}{\sqrt{1 - d_0^2}}
    \cdot \dist_{2 \to \infty}(Q_0, V) \right] \\
                                     & + 
        \left(\frac{\lambda_{r + 2}}{\lambda_r}\right)^t
        \frac{d_0}{\sqrt{1 - d_0^2}},
    \end{aligned}
    \end{align}
    where $\mu$ is the coherence of $V$.
\end{proposition}
Typically, we expect that $\dist_{2\to\infty}(Q_0, V) \ll \dist_{2}(Q_0, V)$,
since otherwise the error is highly localized in just a few rows of the matrix.
Therefore,~\Cref{prop:2inf-convergence-general,prop:2infnorm-convergence-single}
show that we can achieve significant practical improvements in that regime
(recall that convergence analysis with respect to the spectral norm gives
a rate of $\left( \sfrac{\lambda_{r+1}}{\lambda_r} \right)^t \frac{d_0}{\sqrt{1 - d_0^2}}$
\cite{GVL13}).
\Cref{sec:applications} illustrates this concept in practical examples.

\section{Stopping criteria}
\label{sec:stopping}
In this section, we propose and analyze a stopping criterion for tracking convergence
with respect to the $2 \to \infty$ norm. Notably, this
stopping criterion is generic and applicable to any iterative method for computing an invariant subspace.\footnote{This includes \Cref{alg:subspace-iteration} and other common methods such as (block) Lanczos.} Suppose that we
have
\[
    A Q - Q S = E, \quad \norm{E}_2 \leq \varepsilon, \quad
    Q \in \Obb_{n, r}, \quad S = S^\T.
\]
Then it is well-known~\citep[Theorem 8.1.13]{GVL13} that there exist
$\mu_1, \dots, \mu_{r} \in \Lambda(A)$ such that
\(
    \abs{\mu_k - \lambda_k(S)} \leq \sqrt{2} \varepsilon, \quad
    \forall k \in [r].
\)
This provides a handy criterion for testing convergence of eigenvalues, by
setting $S = D_{t}$, the diagonal matrix of {approximate eigenvalues} at the
$t^{\text{th}}$ step and $Q = Q_t$, the orthogonal matrix of {approximate eigenvectors}.
The following lemma is straightforward to show.
\begin{lemma}
    \label{lemma:invariant-subspace-sylvester}
    Suppose that $A = A^\T \in \Rbb^{n \times n}$ satisfies
    \( A Q - Q S = E, \; Q \in \Obb_{n, r} \),
    for some \textbf{diagonal} matrix $S$. Then $Q$ is an invariant
    subspace of the matrix $A - E Q^\T$.
\end{lemma}
We demonstrate that checking $\norm{A Q - Q S}$ leads to an
appropriate stopping criterion for iterative methods, and simplifies under standard
incoherence assumptions. The proof of Proposition~\ref{prop:stopping-criterion}
crucially relies on a perturbation bound from~\cite{DamleSun19}.\footnote{As the perturbed matrix is non-normal, an eigengap
condition does not suffice to guarantee that $V$ is the leading invariant
subspace of the perturbed matrix. To invoke
Proposition~\ref{prop:stopping-criterion} with the approximate eigenvectors in
the place of $Q$, one relies on the fact that $Q$ approaches the leading
eigenvector matrix $V$ by
convergence theory of subspace iteration. For more details, we refer the reader
to~\cref{appendix:eigenvalue-localization}.}
\begin{proposition}
	\label{prop:stopping-criterion}
    Assume that $A$ is symmetric with $V$ as its dominant subspace and $V_{\perp}$
	spans the orthogonal complement of $V$, with $V \in \Obb_{n, r}$; furthermore, suppose
	that $A$ satisfies the conditions of
	Lemma~\ref{lemma:invariant-subspace-sylvester} for some $Q$ and let
    $\gap := \min\set{\lambda_r(A) - \lambda_{r+1}(A), \sep_{(2, \infty), V_{\perp}}(\Lambda,
    V_{\perp} \Lambda_{\perp} V_{\perp}^\T)}$. Then, if
	$Q$ is the \textit{leading} invariant subspace of $A - EQ^\T$ and
	$\norm{E}_2 \leq \frac{\gap}{5}$, we have
	\begin{align*}
		\dist_{2 \to \infty}(V, Q) \leq
			8 \norm{V}_{2 \to \infty} \left(
			\frac{\norm{E}_2}{\lambda_r - \lambda_{r+1}} \right)^2
            + 2 \norm{V_{\perp} V_{\perp}^\T}_{\infty}
			\frac{\norm{E}_{2 \to \infty}}{\gap} \cdot
			\left(1 + \frac{2 \norm{E}_2}{\lambda_r - \lambda_{r + 1}} \right).
	\end{align*}
\end{proposition}
\begin{corollary}
	\label{corl:stopping-incoherence}
    Suppose that $V \in \Obb_{n, r}$ with coherence $\mu$ and that the conditions
    of~\Cref{lemma:invariant-subspace-sylvester} are satisfied with
	$\norm{E}_2 \leq \varepsilon_1, \; \norm{E}_{2 \to \infty} \leq \varepsilon_2$.
	Then the approximate eigenvector matrix $Q$ satisfies
	\begin{equation}
        \dist_{2 \to \infty}(V, Q) \leq 8 \mu \sqrt{\frac{r}{n}}
            \left(\frac{\varepsilon_1}{\lambda_r - \lambda_{r + 1}
            }\right)^2
            + 2 \frac{1 + \mu \sqrt{r}}{\gap} \cdot
            \left( \varepsilon_2 + 2\frac{\varepsilon_1 \varepsilon_2}{
                \lambda_r - \lambda_{r + 1}}
            \right),
            \label{eq:stopping-incoherence}
    \end{equation}
    with $\gap$ defined as in~\Cref{prop:stopping-criterion}.
\end{corollary}

\paragraph{Practical issues.}
Checking the criterion of~\Cref{corl:stopping-incoherence} requires
computing $\norm{E}_{2}$, $\norm{E}_{2 \to \infty}$ and estimating $\gap$.
The first two terms are straightforward.
To estimate $\gap$ in practice,
we assume that $\sep_{2 \to \infty, V_{\perp}}(\Lambda, V_{\perp} \Lambda_{\perp} V_{\perp}^\T)$ is
a small multiple of the $\lambda_r - \lambda_{r + 1}$, motivated by
the observation that $\sep_{2 \to \infty}$ is \textit{at worst} a factor of
$\frac{1}{\sqrt{n}}$ smaller than the eigengap~\cite[Lemma 2.4]{DamleSun19};
moreover, this $\frac{1}{\sqrt{n}}$ factor is typically loose.
To estimate $\lambda_r - \lambda_{r + 1}$, we may use a combination of techniques, such as
augmenting the ``seed'' subspace by a constant number of columns and
setting $\abs{\lambda_r - \lambda_{r+1}} \approx{\hat{\lambda}_r -
\hat{\lambda}_{r+1}}$ -- where $\hat{\lambda}_i = \lambda_i(Q^{\T} A Q)$ are
the approximate eigenvalues -- as it is well known that eigenvalue estimates
converge at a quadratic rate for symmetric matrices~\cite{Stewart1969}.

In the absence of incoherence information, it is not possible to evaluate~\cref{eq:stopping-incoherence}, and
we may instead  replace all quantities in the residual by estimates (which is common
practice for unknown quantities in standard eigensolvers).
For any $B$, $\norm{BQ_t}_{2 \to \infty} \approx
\norm{BV}_{2 \to \infty}$ (by~\cite[Prop. 6.6]{CapeTangPriebe19} and since $Q_t Q_t^\T
\approx VV^\T$ after sufficient progress).
Similar arguments for the other terms yield an approximated residual:

\begin{align}
    \label{eq:crit-no-incoherence}
    \begin{aligned}
    \mathsf{res}_{2 \to \infty}(t) := 8 \norm{Q_t}_{2 \to \infty} \left(
        \frac{\norm{E}_2}{\lambda_r - \lambda_{r + 1}}\right)^2 +
    \frac{2 \norm{(I - Q_t Q_t^\T) E}_{2 \to \infty}}{\gap} \cdot
    \left(1 + \frac{2 \norm{E}_2}{\lambda_r - \lambda_{r+1}}\right).
    \end{aligned}
\end{align}
The main drawback of using~\Cref{eq:crit-no-incoherence} is that the
substitutions used above are not accurate until $Q_t Q_t^\T$ is sufficiently
close to $V V^\T$. This is observed empirically in~\Cref{sec:applications},
as $\mathsf{res}_{2 \to \infty}(t)$ is looser than average in the first few iterations.

Another practical concern is evaluating the quality of the bound in
\cref{corl:stopping-incoherence}; 
there is no known method for computing the $2\to\infty$ subspace distance
$\min_{Z \in \Obb_{r}} \norm{\hat{V} - V Z}_{2 \to \infty}$
in closed form or via some globally convergent iterative method. However,
rather than computing $Z_{\star} = \argmin_{Z \in \Obb_r} \norm{\hat{V} - V Z}_{2\to\infty}$,
we can instead substitute 
\(
    Z_{F} = \argmin_{Z \in \Obb_{r}} \norm{\hat{V} - V Z}_F,
\)
the minimizer of the so-called \textit{orthogonal
Procrustes problem}, whose solution can be obtained via the SVD of $V^\T
\hat{V}$~\cite{Higham88}, as a proxy for tracking the behavior of the
$\twoinf$ distance; this is precisely the solution used by~\cite{DamleSun19}
to study perturbations on the $\twoinf$ distance.
Via standard arguments, we are able to show that the aforementioned
proxy $\norm{\hat{V} - V Z_F}_{2 \to \infty}$ enjoys a similar convergence guarantee
with an additional multiplicative factor of $\sqrt{r}$, which is typically
negligible compared to $n$ -- the details are in the supplementary material.


\section{Applications}
\label{sec:applications}
In this section, we present a set of numerical experiments illustrating the
results of our analysis in practice, as well as the advantages of the proposed
stopping criterion. Importantly, in our applications, \emph{entry-wise} error
is the natural criterion, often because what matters for the downstream
task is an ordering induced by computed eigenvectors. The supplementary material
contains more details about the implementation and the experimental setup.

\paragraph{Synthetic examples.}
To verify our theory and get a sense of the tightness of our bounds on convergence rates,
we first test on synthetic data. To this end, we generate matrices as follows, given a
pair of matrix and subspace dimensions $(n, r)$:
\begin{enumerate}
    \item Sample a matrix from $\Obb_{n}$ uniformly at random
    (see~\cite{Mezzadri07} for details) and select $r$ of its columns uniformly at
    random to form $V$.
    \item generate $\lambda_i \equiv \rho^{i - 1}$, for a decay factor $\rho = 0.95$.
    \item Form $A = \bmx{V & V_{\perp}} \Lambda \bmx{V & V_{\perp}}^\T$, where
    $V_{\perp}$ is initialized as a random subset of the columns of the
    identity matrix, and subsequently orthogonalized against $V$.
\end{enumerate}

\begin{figure*}[tb]
    \centering
    \begin{minipage}{0.33\textwidth}
        \includegraphics[width=\linewidth]{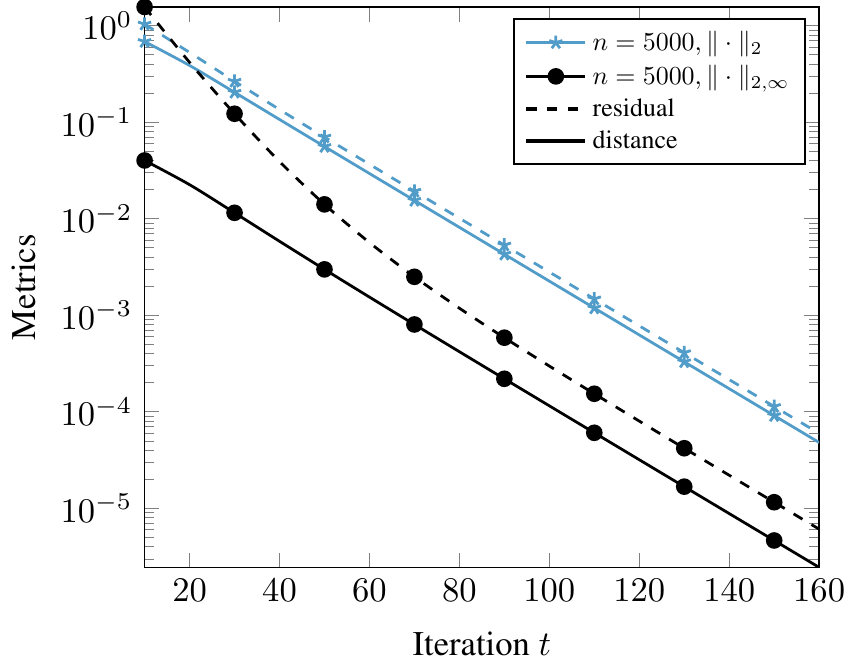}
    \end{minipage}~
    \begin{minipage}{0.33\textwidth}
        \includegraphics[width=\linewidth]{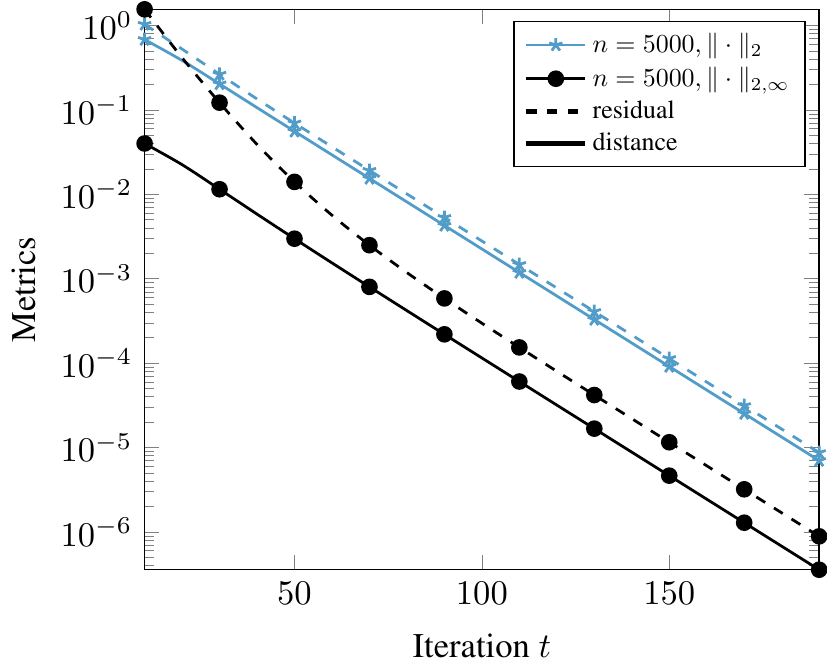}
    \end{minipage}~
    \begin{minipage}{0.33\textwidth}
        \includegraphics[width=\linewidth]{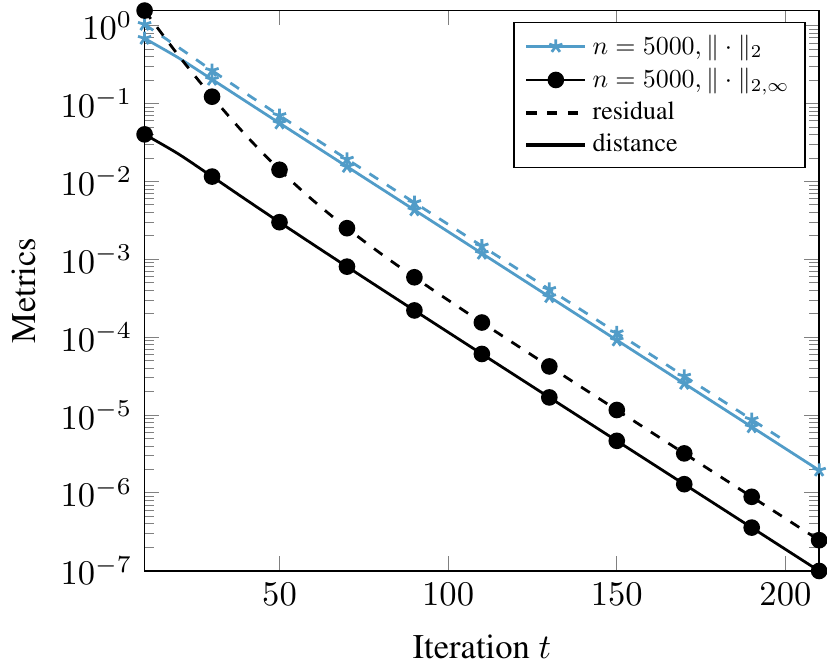}
    \end{minipage}
    \caption{Distances (\textbf{solid} lines) and residuals (\textbf{dashed}
    lines) for synthetic examples with $r = 50$
    and target accuracies $\varepsilon = 10^{-4}$ (\textbf{left}),
    $\varepsilon = 10^{-5}$ (\textbf{middle}) and
    $\varepsilon = 10^{-6}$ (\textbf{right}). Each plot corresponds to an
independently generated synthetic example.}
    \label{fig:synthetic}
\end{figure*}

We compare distances and residuals for synthetic examples with $n = 5000$ and
$r = 50$ and various stopping thresholds $\varepsilon$ for the residuals
(\Cref{fig:synthetic}). Each plot in \cref{fig:synthetic} corresponds to
a different matrix generated independently according to the aforementioned scheme.
While the $\ell_2$ norm residual closely tracks the corresponding distance, the
residual from~\Cref{eq:crit-no-incoherence} overshoots by a small multiplicative
factor, suggesting that the large constants in \Cref{prop:stopping-criterion}
may only be necessary in pathological cases and could be reduced in practice.
Moreover, the $\twoinf$ norm residual can substantially overestimate the actual
distance in the first few iterations, as the estimate of
\Cref{eq:crit-no-incoherence} depends on $Q_t Q_t^\T$ not being ``too far''
from $V V^\T$.  The gap narrows after a few dozen iterations.

\begin{figure*}[tb]
    \centering
    \begin{minipage}{0.33\textwidth}
        \includegraphics[width=\linewidth]{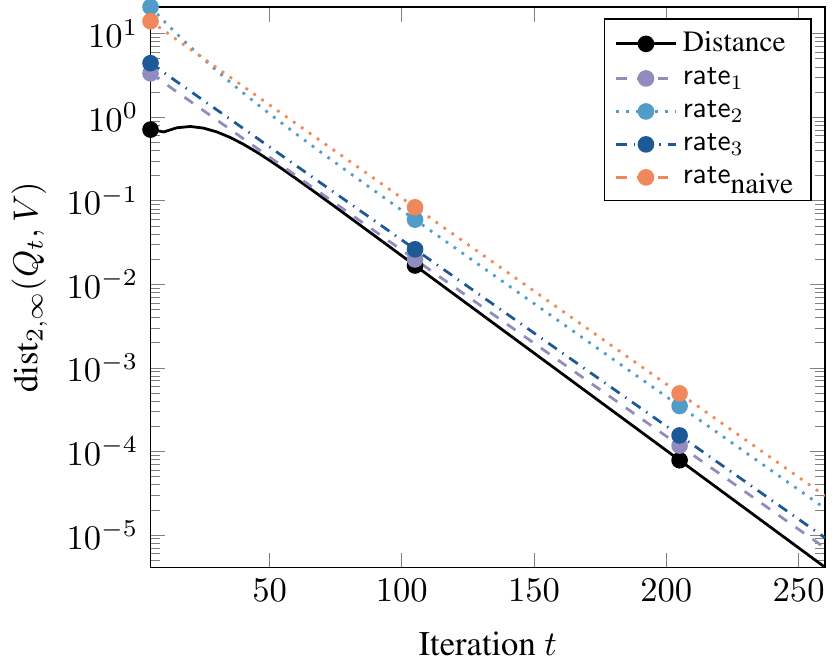}
    \end{minipage}~
    \begin{minipage}{0.33\textwidth}
        \includegraphics[width=\linewidth]{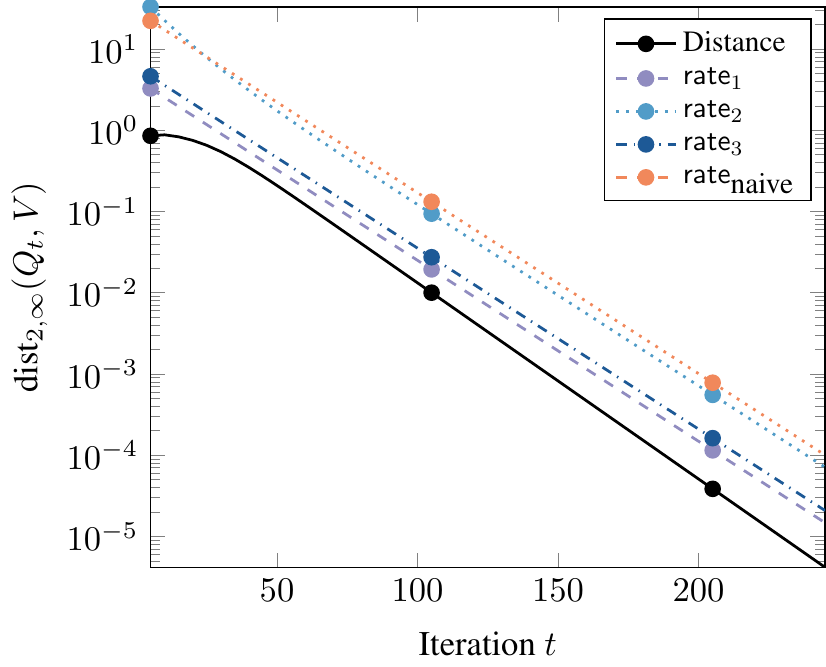}
    \end{minipage}~
    \begin{minipage}{0.33\textwidth}
        \includegraphics[width=\linewidth]{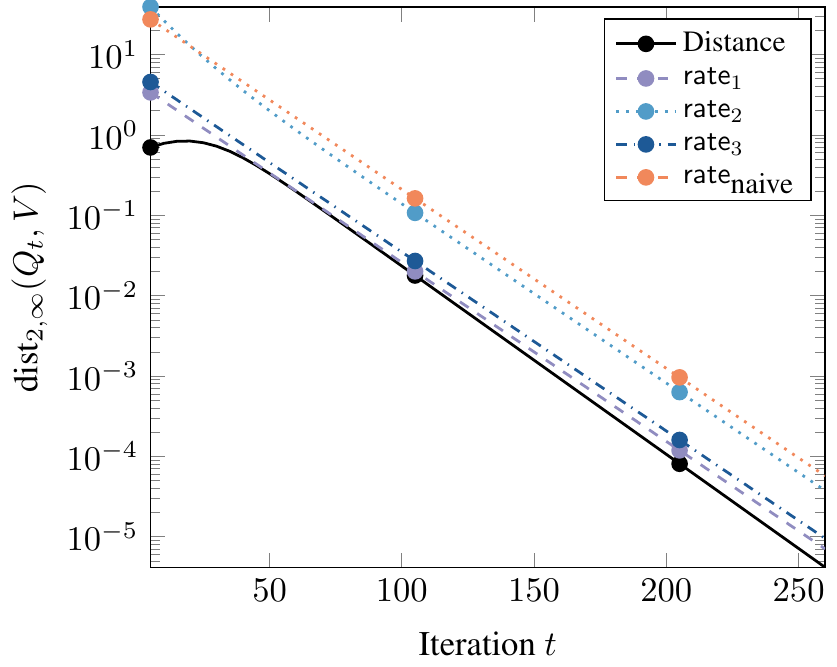}
    \end{minipage}
\caption{Distance (\textbf{solid} lines) and convergence rates
  from \Cref{eq:rates} for matrix and subspace dimensions $(n, r) = (1000, 10)$
  (\textbf{left});
  $(3500, 15)$ (\textbf{middle}); and $(8000, 20)$ (\textbf{right}).
  Our $\mathsf{rate}_3$ from \cref{prop:2inf-convergence-general} tracks the
  ``idealized'' rate $\mathsf{rate}_1$ closely in the synthetic data examples.
  }
\label{fig:convergence-rates}
\end{figure*}

In addition, we examine the looseness of the bounds from
\Cref{prop:2inf-convergence-general,prop:2infnorm-convergence-single} for the
same experiment (\Cref{fig:convergence-rates}).  We evaluate the following
rates:
\begin{align}
    \begin{aligned}
    \mathsf{rate}_1(t) &:= \left( \frac{\lambda_{r+1}}{\lambda_r} \right)^t
    \cdot \frac{\dist_{2 \to \infty}(Q_0, V)}{\sqrt{1 - d_0^2}}, \quad
    \mathsf{rate}_2(t) := \text{ rate
    from~\Cref{prop:2infnorm-convergence-single}}, \\
    \mathsf{rate}_3(t) &:= \text{ rate
    from~\Cref{prop:2inf-convergence-general}}, \quad
    \mathsf{rate}_{\textrm{naive}}(t) := \left(\frac{\lambda_{r + 1}}{\lambda_r}\right)^t
    \frac{d_0}{\sqrt{1 - d_0^2}}
    \end{aligned}
    \label{eq:rates}
\end{align}
Here, $\mathsf{rate}_1$ is an idealized rate that mirrors classical convergence
results for the $\ell_2$ norm \citep[Theorem 8.2.2]{GVL13}; on the other hand, the
naive rate just measures the $\ell_2$ subspace distance. In all the synthetic
examples we generated, Assumption~\ref{asm:norm-bound} was verified to hold
with constant $C < 2$ for all elapsed iterations $t$.

Remarkably, for a range of dimensions $n$ and $r$ we find that $\mathsf{rate}_3$
(which uses \cref{prop:2inf-convergence-general}) closely tracks the
``idealized'' $\mathsf{rate}_1$ on these synthetic matrices
(\cref{fig:convergence-rates}). Also, $\mathsf{rate}_2$ (which uses
\Cref{prop:2infnorm-convergence-single}) is a looser upper bound. This agrees
with our theoretical analysis, as $\lambda_{r+2}$ is only moderately smaller
than $\lambda_{r+1}$ in our synthetic matrix construction.  Finally, as
expected, the naive rate is the loosest bound.

\paragraph{Eigenvector centrality.}
Next, we develop an experiment for network centrality, where the
task is to measure the influence of nodes in a graph~\cite{Newman08}.
Each node is assigned a score, which is a function of the graph topology,
and a typical underlying assumption is that a node with a high
score contributes a larger influence to its adjacent ones.
Here, we consider \emph{eigenvector centrality}, which is one the standard
measures in network science.
Given a graph $G = (V, E)$; the eigenvector centrality score of a node
$u$, $x_u > 0$, is defined as a solution to the following equation:
\begin{equation}
    x_u := \frac{1}{\lambda} \sum_{v \in V} A_{uv} x_v, \;
    A_{uv} := \begin{cases}
        1, & \text{if $u$ links to $v$} \\
        0, & \text{otherwise}
    \end{cases},
    \label{eq:rel-eig-centrality}
\end{equation}
where $\lambda$ is a proportionality constant. Here, node $u$'s scores depend
linearly on all of its neighbors' scores.
Under the positivity requirement of $x_u$ and provided that the graph is
connected and non-bipartite, 
rearranging and the Perron-Frobenius theorem show that $x = v_1$, the leading eigenvector of $A$ (up to scaling).
To determine the most influential nodes, we are typically interested in the
induced \emph{ordering} of nodes and not the actual scores themselves.
Therefore, the $\ell_{2 \to \infty}$ distance, which measures $\norm{v_1 - \hat{v}_1}_{\infty}$,
is more appropriate than $\norm{v_1 - \hat{v}_1}_2$ as a proxy for the quality of the estimate $\hat{v}_1$.
To get a correct ranking result, it suffices to have
$\norm{v_1 - \hat{v}_1}_{\infty} < (1/2) \cdot \min_{i, j} \abs{v_i - v_j}$.
On the other hand, $\norm{\hat{v}_1 - v_1}_2$ does not have an
interpretable criterion.

We demonstrate the above principle by comparing two stopping
criteria: the criterion from~\Cref{eq:crit-no-incoherence} with a specified
threshold $\varepsilon$ against the ``naive'' way of stopping when $\norm{A
\hat{v}_1 - \hat{\lambda} \hat{v}_1}_2 \leq \hat{\lambda} \varepsilon$,
where $\hat{\lambda}$ is the current eigenvalue estimate, via the two following
stopping times:
\begin{align}
    \begin{aligned}
    t_{\textrm{comp}} &:= \min\set{t > 0 \mmid
    \mathsf{res}_{2 \to \infty}(t) \leq \varepsilon} \\
    t_{\textrm{naive}} &:= \min\set{t > 0 \mmid
    \norm{A \hat{V}_{:, j} - \hat{\lambda}_j \hat{V}_{:, j}}
\leq \varepsilon \hat{\lambda}_j, \forall j \in \set{1, \dots, r}}
    \end{aligned}.
    \label{eq:tcomp-tnaive}
\end{align}
For a user-specified tolerance $\varepsilon$, we
expect that using our $\twoinf$ error measurements and our corresponding
stopping criteria will tell us that we can be confident in our solution
much more quickly.
\begin{table}[tb]
    \centering
    \caption{Summary statistics of network datasets.}
    \begin{tabular}{l c c c} \toprule
        \textbf{Dataset}  & \textbf{Citation} & \textbf{\# nodes}    & \textbf{\# edges} \\ \midrule
        \textsc{ca-HepPh} & \multirow{2}{*}{\cite{leskovec2005graphs}} &  11204 & 117649 \\
        \textsc{ca-AstroPh} & & 17903 & 197031 \\ \midrule
        \textsc{gemsec-facebook-artist} & \cite{RDDS19} & 50515 & 819306 \\ \midrule
        \textsc{com-DBLP} & \multirow{2}{*}{\cite{YangLesk12}} & 317080            & 1049866 \\
            \textsc{com-LiveJournal} & & 3997962           & 34681189 \\ \bottomrule
    \end{tabular}
    \label{tab:eigvec-centrality-data}
\end{table}
\begin{figure*}[tb]
  \centering
  \begin{minipage}{0.44\linewidth}
      \includegraphics[width=\textwidth]{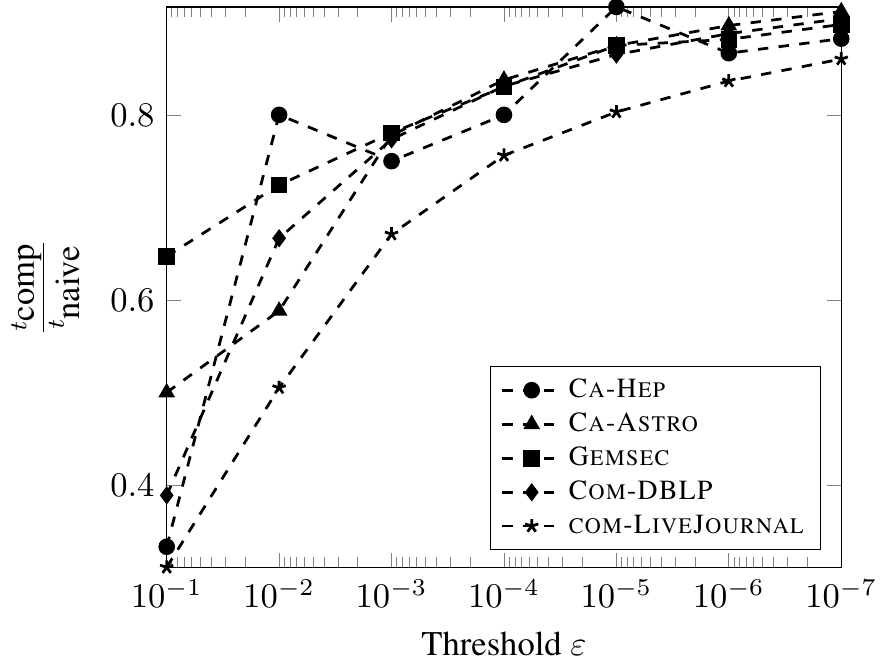}
  \end{minipage}\quad
  \begin{minipage}{0.44\linewidth}
      \includegraphics[width=\textwidth]{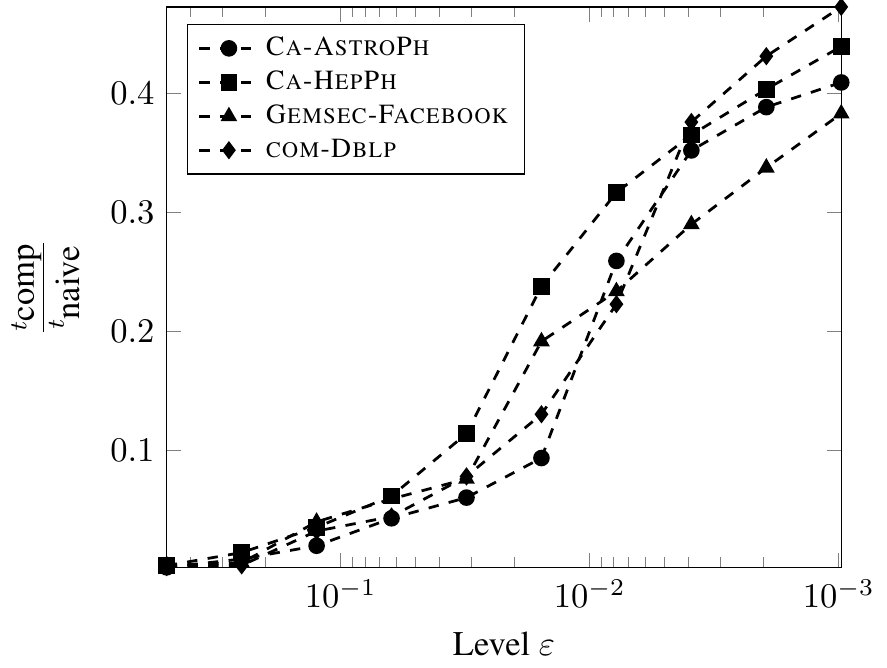}
  \end{minipage}
  \caption{Ratio of the number of iterations needed to satisfy the two stopping
      criteria outlined in~\eqref{eq:tcomp-tnaive}, for thresholds $\varepsilon = 10^{-k}$, for
  computing eigenvector centrality to find the $\floor{\sqrt{n}}$ most influential
  nodes (\textbf{left}) and computing the leading $r$ eigenvectors for spectral clustering (\textbf{right}). Our analysis and stopping criteria enable significantly fewer iterations.}
  \label{fig:evec-ratios}
\end{figure*}
%
This is indeed the case --- using our methodology provides a substantial
reduction in computation time on a variety of real-world graphs, whose
summary statistics are in \cref{tab:eigvec-centrality-data}.
\Cref{fig:evec-ratios} (left) shows the ratio between the two
quantities $t_{\textrm{comp}}$ and $t_{\textrm{naive}}$, defined as in~\cref{eq:tcomp-tnaive}.
In the low-to-medium accuracy regimes, using our
stopping method results in \textbf{at least a 20--40\% reduction in the number of
iterations needed}. In this regime, the ranking induced by the approximate
eigenvector had typically already converged to the ``true'' ordering obtained
by computing the eigenvector to machine precision.
\paragraph{Spectral clustering in graphs.}
Another downstream task employing invariant subspaces is spectral clustering,
which we study here as a way to partition a graph into well-separated
``communities'' or ``clusters.''  The standard pipeline
is to compute the leading $r$-dimensional eigenspace of the normalized
adjacency matrix, where $r$ is the desired number of clusters,
The resulting eigenvector matrix provides an $r$-dimensional embedding for each node,
which is subsequently fed to a point cloud clustering algorithm such as
\texttt{k-means}~\cite{VonLux07}. For our experiment, we use the
\textit{deterministic} QR-based algorithm from~\cite{DMY18} on the same set of
real-world graphs that we used for eigenvector centrality.

In this setup, the eigenvectors (more carefully, a rotation of them) are
approximate cluster indicators. Indeed, spectral clustering on graphs is often
derived from a continuous relaxation of a combinatorial objective based on these
indicators~\cite{VonLux07}. Thus, we are once again interpreting the
eigenvectors entry-wise, and $\twoinf$ error is a more sensible metric than $\ell_2$ error,
This fact has been used to analyze random graph models with cluster structure~\cite{lyzinski2014perfect}.

In the same manner as the eigenvector centrality experiment, we compare the
ratio of iteration counts: $t_{\textrm{comp}}$ over $t_{\textrm{naive}}$, as
defined in~\Cref{eq:tcomp-tnaive} (\Cref{fig:evec-ratios}, right).  In this
case, we see even larger savings. For $\varepsilon$ around $10^{-2}$, our
stopping criterion results in 70--80\% savings in computation time.  While this
approximation level may seem crude at first, we can measure the performance of
the algorithms in terms of the normalized cut metric, for which spectral
clustering is a continuous relaxation~\cite{VonLux07}. We find that by the time
we reach residual level $\varepsilon = 10^{-2}$, the cut value found using the
approximate subspace is essentially the same as the one using the subspace
computed to numerical precision. Further details about the experiment are
provided in the supplementary material.

\paragraph{Spectral bipartitioning and sweep cuts.}
Another spectral method for finding clusters in graphs is spectral bipartitioning,
which aims to find a single cluster of nodes $S$ with small conductance
$\phi(S)$:
\begin{equation*}
    \phi(S) := \frac{\sum_{i \in S, j \notin S} A_{ij}}{\min(A(S), A(S^c))}, \;
    A(S) := \sum_{i \in S} \sum_{j \in V} A_{ij}.
\end{equation*}
The conductance objective is a standard measure for identifying a good
cluster of nodes~\cite{schaeffer2007graph,LLDM08}: if $\phi(S)$ is small, there
are not many edges leaving $S$ and there are many edges contained in $S$.

Minimizing $\phi(S)$ is NP-hard, but a spectral method using the eigenvector
$v_2$ corresponding to the second largest eigenvalue of the normalized adjacency
matrix, often called the \textit{Fiedler vector}~\cite{Fiedler73}, provides
guarantees. To find the a set with
small conductance, the method uses the so-called ``sweep cut''. After scaling
$v_2$ by the inverse square root of degrees, we sort the nodes by their value in
the eigenvector, and then consider the top-$k$ nodes as a candidate set $S$ for
all values of $k$. The value of $k$ that gives the smallest conductance produces
a set $S$ satisfying $\phi(S) \le 2\sqrt{\min_{S'} \phi(S')}$, which is the
celebrated Cheeger inequality~\cite{chung1997spectral}.

\begin{figure}[tb]
    \centering
    \begin{minipage}{0.5\textwidth}
      \includegraphics[width=\linewidth]{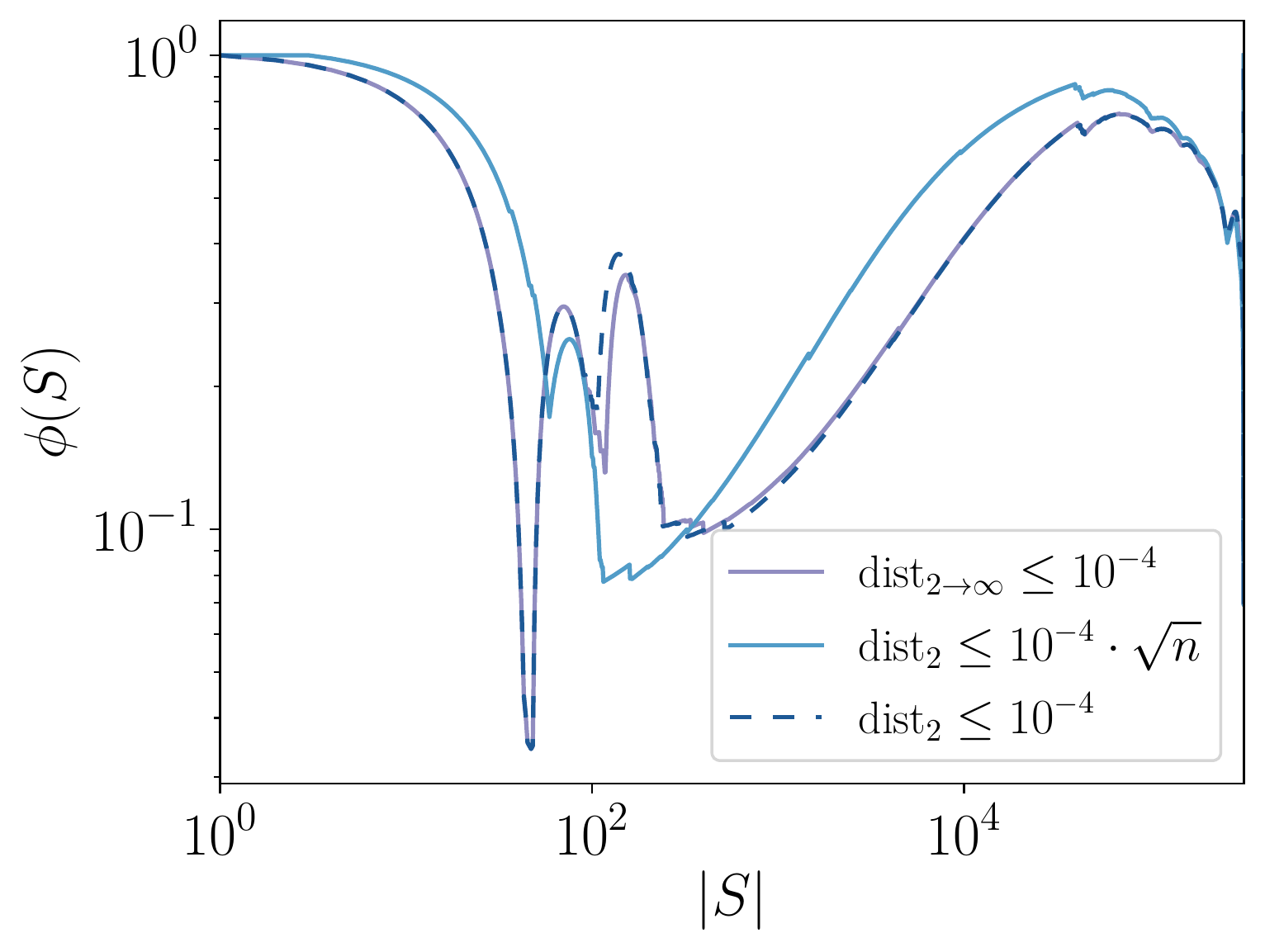}
    \end{minipage}
    \hfill
    \begin{minipage}{0.45\textwidth}
      \caption{Sweep cut profile (cut conductance vs.\ cardinality)
      for \textsc{com-Dblp}. For a fixed $\varepsilon$, our $\twoinf$
      stopping criterion leads to faster convergence. Increasing the tolerance
      for $\ell_2$ by the norm equivalence factor produces lower-quality
      solutions. Here $t_{\textrm{comp}} = 1135$ vs. $t_{\textrm{naive}} = 1378$
      iterations.}
      \label{fig:ncp}
    \end{minipage}
\end{figure}

As in the case of eigenvector centrality, what matters is the \emph{ordering}
induced by the eigenvector, making a $\twoinf$ stopping criterion more appropriate.
As a heuristic, one might consider just making $\ell_2$ tolerance larger (by the
norm equivalence factor) using a level of $\varepsilon$ for the
$\twoinf$ distance and $\varepsilon \cdot \sqrt{n}$ for the $\ell_2$ distance.
However, this can substantially reduce the solution quality.
This is illustrated in \cref{fig:ncp}, where we plot the conductance values
obtained in the sweep cut as a function of the size of the set on
\textsc{com-Dblp}. This is a sweep cut approximation of
a network community profile plot~\cite{benson2016higher,LLDM08},
which visualizes cluster structure at different scales. Using the
naive $\ell_2$ stopping criterion provides the same solution quality but
requires more iterations.
In the case of $\varepsilon = 10^{-4}$ in \cref{fig:ncp}, our methods produce 20\% computational savings.
Finally, the heuristic $\varepsilon \cdot \sqrt{n}$ tolerance for the $\ell_2$ stopping criterion
produces a cruder solution and finds a set with larger conductance.

\section{Conclusions}
\label{sec:conclusions}
The broad applicability of spectral methods, coupled with the prevalence
of entry-wise / row-wise interpretations of eigenspaces strongly motivates
imbuing our computational methods with appropriate stopping criteria.
Our theoretical results demonstrate just how
much smaller the $\norm{\cdot}_{2\to \infty}$ subspace distance can be than
traditional measures, an observation supported by experiment.
In fact, the accuracy with which we compute eigenvectors can have a non-trivial impact on
downstream applications\;\textemdash\; if we would like use fewer iterations to
save time we must do so carefully, and our new stopping criterion provides an
easy to implement way to do this that comes at essentially no cost and with
strong guarantees.

From a theoretical perspective, it may seem sufficient to use norm equivalence
and simply appeal to spectral norm convergence, which can incur an extra $\cO(\log n)$
factor at most when computing subspaces. However, such reasoning only applies
to the very-high-accuracy regime.
As demonstrated by our experiments, moderate
levels of accuracy often suffice for downstream applications, in which case our
stopping criterion allows for highly nontrivial computational savings (up to 70\%
fewer iterations).

\section*{Broader impact}
Due to the pervasiveness of spectral methods in machine learning and data mining,
our results may be embedded in applications having a wide range of ethical and
societal consequences. Indeed, given the fact that eigensolvers are typically
used as linear algebra primitives, our work ``inherits'' the ethical and societal
consequences of the context in which its results are applied, as well as the
potential implications of ``failure'' (\emph{e.g.}, if our stopping criterion
severely underestimates the true approximation error).

\section*{Acknowledgements \& Funding}
We would like to thank the anonymous reviewers for their valuable feedback,
which helped improve the presentation of this work.

This research was supported by NSF Award DMS-1830274, ARO Award W911NF19-1-0057, and ARO MURI.

%
%
%
%
%

\bibliographystyle{plain}
\bibliography{references}

\appendix

\section{Additional experimental details}
This section documents hyperparameters and other design choices used to run the
experiments in Section 4 of the main paper.

\subsection{Eigenvector centrality}
In all eigenvector centrality experiments, we isolate the largest connected
component of the input graph and work exclusively within that component. We
work with the unnormalized version of the adjacency matrix, since the
normalized version admits the vector $d^{1/2} := \diag(\sqrt{d_1}, \dots,
\sqrt{d_n})$, where $d_i$ is the degree of the $i^{\textrm{th}}$ node, as its
principal eigenvector.

We initialize the estimate $\hat{v}_1 := \frac{1}{\sqrt{n}} \bm{1}$, the
normalized all-ones vector. In the absence of incoherence, we use the
expression of~\Cref{eq:crit-no-incoherence} in the main text to evaluate the
$\ell_{2 \to \infty}$ stopping criterion, and set $\mathsf{gap} = \lambda_1 -
\lambda_2$
using the values returned by \texttt{Arpack}, to ensure a fair comparison.
At each step of the iterative method, we multiply with $\pm 1$ accordingly, to
ensure that all entries of the approximate eigenvector are positive.

\paragraph{Ranking distance.} To measure the ``distance'' between the
approximate ranking produced by our eigenvector estimate, we employ Kendall's
$\tau$ criterion~\cite{Kendall48}. In particular, we define
\begin{equation}
	\dist_{\tau}(v_1, \hat{v}_1) := \frac{1 - \tau(v_1, \hat{v}_1)}{2}
	\label{eq:tau-dist}
\end{equation}
to compare the rankings induced by $v_1$ and $\hat{v}_1$. It is easy to verify
that when the rankings are identical, $\dist_{\tau} = 0$, and when the rankings
are the most dissimilar, $\dist_{\tau} = 1$, since $\tau(v_1, \hat{v}_1) \in
[-1, 1]$.

\begin{figure}[h]
    \centering
    \begin{minipage}{0.4 \textwidth}
        \includegraphics[width=\linewidth]{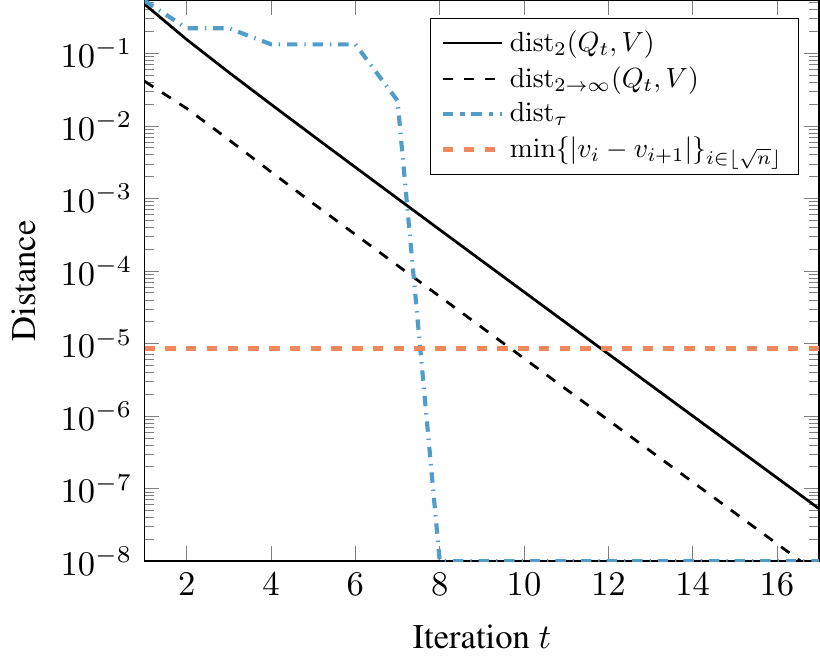}
    \end{minipage} \qquad
    \begin{minipage}{0.4 \textwidth}
        \includegraphics[width=\linewidth]{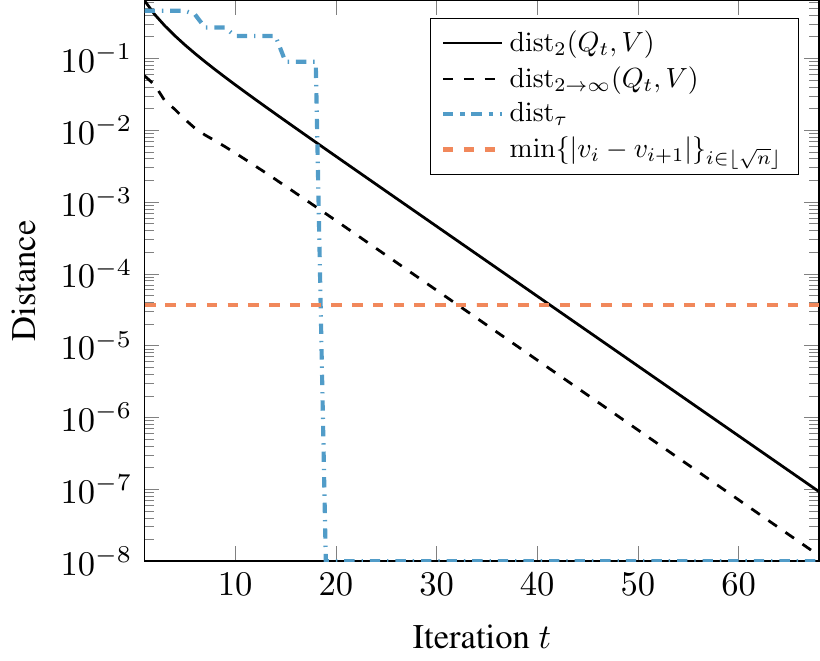}
    \end{minipage}
    \vspace{12pt}
    \begin{minipage}{0.4 \textwidth}
        \includegraphics[width=\linewidth]{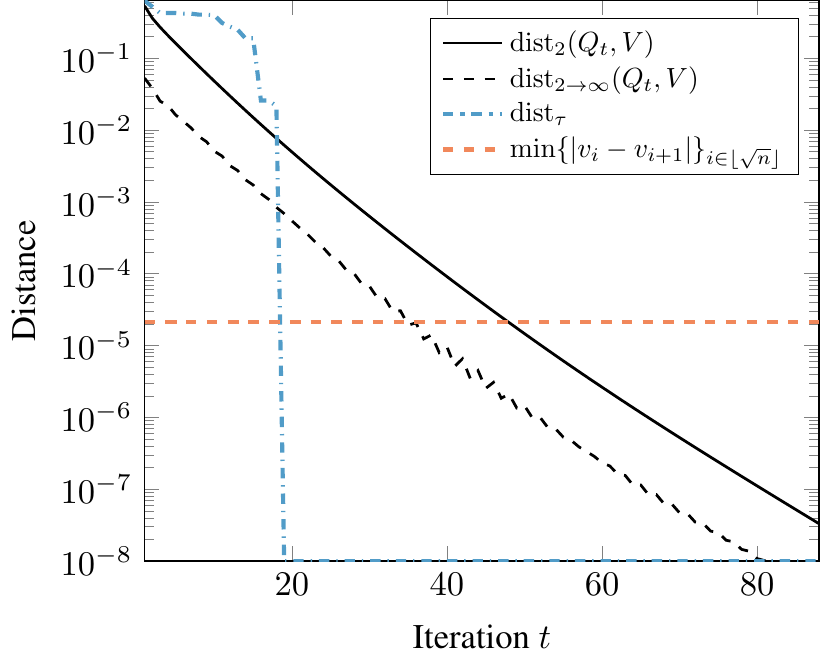}
    \end{minipage} \qquad
    \begin{minipage}{0.4 \textwidth}
        \includegraphics[width=\linewidth]{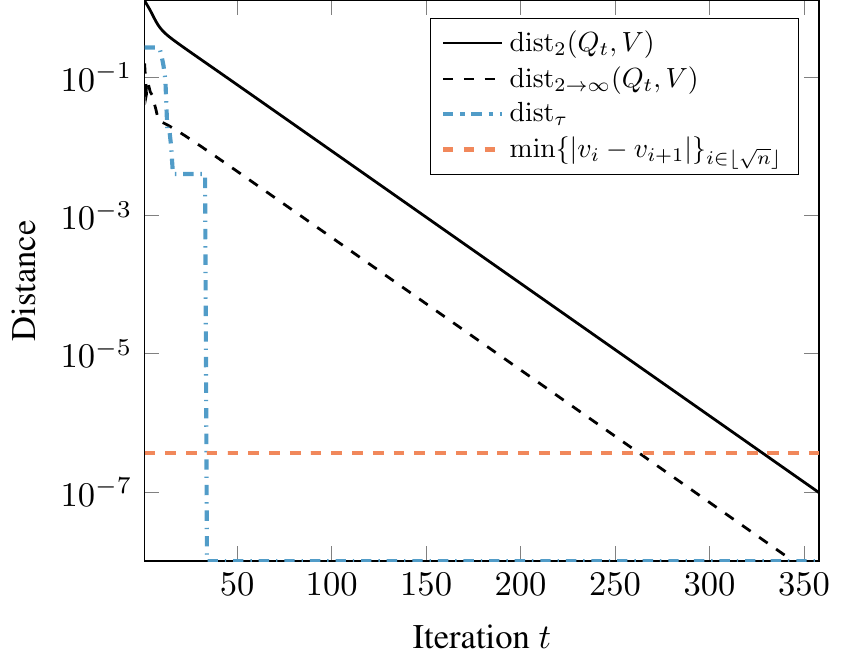}
    \end{minipage}
    \caption{Distance plots for $4$ datasets, for which the top
    $\floor{\sqrt{n}}$ nodes are being ranked. From \textbf{left} to
    \textbf{right}: \textsc{ca-HepPh}, \textsc{ca-AstroPh} (\textbf{top}),
    \textsc{Gemsec}, \textsc{com-LiveJournal} (\textbf{bottom}).}
    \label{fig:tau-plots}
\end{figure}

The convergence plots for $4$ datasets, where we depict the ``oracle'' $\ell_2$
and $\ell_{2 \to \infty}$ subspace distances as well as $\dist_{\tau}$ as a
function of the iteration index $t$, are shown in~\Cref{fig:tau-plots}. In all
cases, we identify the correct ranking when the residual is in the
low-to-moderate accuracy regime ($\varepsilon \leq 10^{-4}$).

\subsection{Spectral clustering}
In this section, we describe the methodology used for the spectral clustering
experiments in the main text. We opt to use the
Algorithm of~\cite{DMY18} which is based
on the column-pivoted QR decomposition of an appropriately defined matrix.
For completeness, the full algorithm is listed in~\Cref{alg:cpqr-clustering}.
Since the algorithm is deterministic, we do not have to worry about randomness
pertaining to initialization (e.g. as in \texttt{kmeans++}), and only run the
experiment once for each configuration of parameters.

\begin{algorithm}
    \caption{CPQR-based clustering}
    \begin{algorithmic}[1]
    \State \textbf{Input}: invariant subspace $V_k \in \Rbb^{n \times r}$
    \State Compute the CPQR factorization
    \[
        V_k^\top \Pi = QR,
    \]
    where $\Pi$ is a column selection matrix.
    \State Let $\cC$ denote the first $k$ columns identified by $\Pi$.
    \State Compute the polar factorization
    \[
        (V_k^\top)_{:, \cC} = U H.
    \]
    \For {$j \in [n]$}
    \State assign node $j$ to cluster
    \[
        C_j := \argmax_{i} \abs{(U V_k^\top)_{i,j}}
    \]
    \EndFor
    \end{algorithmic}
    \label{alg:cpqr-clustering}
\end{algorithm}

\begin{figure}[h]
    \centering
    \begin{minipage}{0.47\textwidth}
    \includegraphics[width=\linewidth]{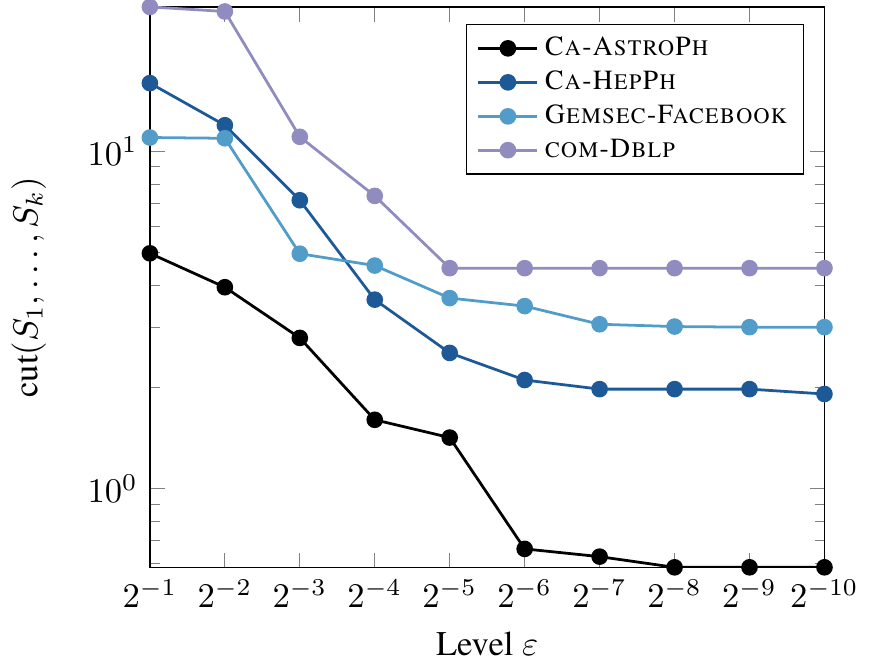}
    \caption{Value of $\textrm{ncut}(S_1, \dots, S_k)$ for various datasets,
    with
    $\hat{V}_k$ computed using subspace iteration until the residual
    drops below level $\varepsilon$, for different values of $\varepsilon$. In
    all cases, the metric stabilizes while in the low accuracy regime ($\varepsilon \approx
    10^{-2}$). Dashed lines indicate the value of $\textrm{ncut}(S_1, \dots,
    S_k)$ found by computing the subspace to machine accuracy.}
    \label{fig:cluster-perf}
    \end{minipage}\qquad
    \begin{minipage}{0.47\textwidth}
    \includegraphics[width=\linewidth]{plots/clusterRatio.pdf}
    \caption{Ratio of iterations required to satisfy $\mathsf{res}_{2 \to
    \infty}(t) \leq \varepsilon$ ($t_{\textrm{comp}}$) over number of iterations
    required to satisfy $\mathsf{res}_2(t) \leq \varepsilon$
    ($t_{\textrm{naive}}$) in the spectral clustering setting, showing
    computational gains of over $50\%$.}
    \end{minipage}
\end{figure}

For all the datasets involved, we hand-pick the target number of clusters $r$
by inspecting the successive ratios of the leading few eigenvalues and setting
$r$ so that the ratio $\frac{\lambda_{r+1}}{\lambda_r}$ is small, but also
taking into account the fact that we don't want $r$ to be too small.
Additionally, we use the regularized version of the normalized adjacency matrix
$A_{\rho}$~\cite{ACBL13}, which augments the adjacency and degree matrices $A,
D$ using a regularization parameter $\rho$:
\begin{equation}
    A_{\rho} := A + \frac{\rho}{n} \bm{1} \bm{1}^\top, \quad
    D_{\rho} := D + \rho
    \label{eq:regadj}
\end{equation}
Following standard practice~\cite{QinRohe13,ZhangRohe18}, we set $\rho$ equal
to a constant which is near the average degree of the graph and then
perform the eigendecomposition of
\[
    \tilde{A}_{\rho} = D_{\rho}^{-1/2} A_{\rho} D_{\rho}^{-1/2} + I,
\]
shifting by $+I$ to ensure that the algebraically largest eigenvalues are
also the largest in magnitude, in order for subspace iteration to be applicable. We summarize the hyperparameter choices for each dataset in~\Cref{tab:spectral-clustering}.
\begin{table}
    \centering
    \caption{Parameters for spectral clustering}
    \begin{tabular}{c | c c}
        \textbf{Dataset} & $r$ & $\tau$ \\ \hline
        \textsc{ca-HepPh} & $17$ & $1.0$ \\
        \textsc{ca-AstroPh} & 6  & $1.0$ \\
        \textsc{GemSec} & 12 & 1.0 \\
        \textsc{Dblp} & 28 & 5.0 \\ \hline
    \end{tabular}
    \label{tab:spectral-clustering}
\end{table}
\begin{figure}[h]
    \centering
    \includegraphics{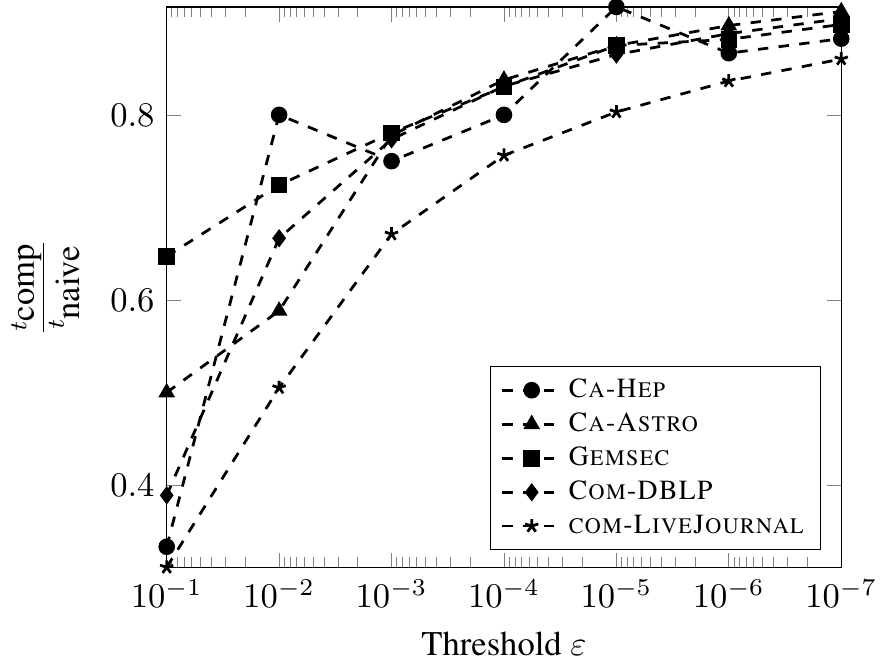}
    \caption{Ratio of iterations required to satisfy $\mathsf{res}_{2 \to
    \infty}(t) \leq \varepsilon$ ($t_{\textrm{comp}}$) over number of iterations
    required to satisfy $\mathsf{res}_2(t) \leq \varepsilon$
    ($t_{\textrm{naive}}$) in eigenvector centrality computations.}
\end{figure}
To evaluate the quality of a given clustering assignment, we use the
\textit{normalized cut} metric. Specifically, given a vertex set $V$ and a
\textit{partition} $(S, S^c)$ such that $V = S \cup S^c$, we define the
conductance of the cut induced by $S$ as
\begin{equation}
    \phi(S) := \frac{\sum_{i \in S, j \notin S} A_{ij}}{A(S)}
    \quad A(S) := \sum_{i \in S} \sum_{j \in V} A_{ij}
    \label{eq:cut}
\end{equation}
Note that in~\eqref{eq:cut}, $A$ refers to the \textbf{unnormalized} adjacency
matrix, with $A_{ij} = A_{ji} = 1$ if the edge $(i, j)$ exists in the graph, and
$0$ otherwise.
Then any clustering assignment with $k$ clusters induces $k$ partitions
$\set{(S_k, S_k^c)}$, for which the normalized cut metric is defined as
\begin{equation}
    \textrm{ncut}(S_1, \dots, S_k) := \frac{1}{2} \sum_{i = 1}^k \phi(S_k).
\end{equation}
\Cref{fig:cluster-perf} depicts the value of $\textrm{ncut}(S_1, \dots, S_k)$
when the input to~\Cref{alg:cpqr-clustering} is computed using subspace
iteration, using the proposed stopping criterion, for different levels
$\varepsilon$.
Having established that low-to-moderate accuracy is sufficient for this problem,
we plot the ratio of $t_{\textrm{comp}}$ over $t_{\textrm{naive}}$; the former
is the number of iterations required to satisfy $\mathsf{res}_{2 \to \infty}(t)
\leq \varepsilon$, while the latter is the number of iterations required to
satisfy $\mathsf{res}_2(t) := \norm{A \hat{v}_t - \hat{\lambda}_t \hat{v}_t}_2
\leq \hat{\lambda}_t \varepsilon$. We observe computational gains of over
$50\%$ in all cases.

\subsection{Empirically verifying Assumption~\ref{asm:norm-bound}}
We verified that Assumption~\ref{asm:norm-bound} from the main text holds in practice for the real
world datasets used in the experimental section. Recall that the assumption asks
that the matrix $A = V \Lambda V^\T + V_{\perp} \Lambda_{\perp} V_{\perp}^\T$ 
satisfies:
\begin{equation}
    \norm{V_{\perp} \Lambda_{\perp}^t V_{\perp}^\T}_{\infty}
    \leq C \cdot \norm{\Lambda_{\perp}^t}_2 \cdot \norm{V_{\perp} V_{\perp}^\T}_{\infty},
    \label{eq:asm-2}
\end{equation}
where $C$ is a constant independent of $n$, for all $t \in \Nbb$. First, observe
that for our purposes, we only want this assumption to hold for all $t$ until our
iterative algorithm stops. Since all our experiments take fewer than $T = 1500$
iterations to run, we opt to verify~\eqref{eq:asm-2} for $t \in \set{1, \dots, T}$.
We first rephrase the assumption as
\begin{equation}
    \norm{A^t - V \Lambda^t V^\T}_{\infty}
    \leq C \cdot \lambda^t_{\max}(\Lambda_{\perp}) \cdot \norm{I - V V^\T}_{\infty},
    \label{eq:asm-2-rephrased}
\end{equation}
which can be checked after computing the top $r+1$ eigenvectors and eigenvalues
of $A$; these were computed to machine precision using \texttt{eigs}. For $t = 1$
up to $t = T$, we checked~\eqref{eq:asm-2-rephrased} exhaustively, and output
\[
    C := \sup_{t \in \set{1, \dots, T}} \set{
        \frac{\norm{A^t - V \Lambda^t V^\T}_{\infty}}{
        \lambda_{r+1}^t \norm{I - V V^\T}_{\infty}}}
\]
In all cases, we end up with a constant $C < 2$.

\section{Auxiliary results}

\begin{lemma}[Incoherence]
    \label{lemma:incoherence-compl}
    Consider a subspace $\cV$ of dimension $r$ and a matrix $V \in \Obb_{n,r}$
    whose columns span $\cV$. If $\mu$ is the coherence of $V$, i.e.
    $\norm{V}_{2 \to \infty} \leq \mu \sqrt{\frac{r}{n}}$, then for its
    complementary subspace $\cV_{\perp}$ it holds that
    \[
        \norm{V_{\perp} {V_{\perp}}^\T}_{\infty} \leq (1 + \mu \sqrt{r}).
    \]
\end{lemma}
\begin{proof}
    Observe that $\norm{A}_{\infty} \leq \sqrt{n} \norm{A}_{2 \to
    \infty}$, thus
	\begin{align*}
        \norm{V_{\perp} {V_{\perp}}^\T}_{\infty} &=
        \norm{I - V V^\T}_{\infty} \leq 1 + \norm{V V^\T}_{\infty}
        \leq 1 + \sqrt{n} \norm{VV^\T}_{2 \to \infty}
        \leq 1 + \sqrt{n} \mu \sqrt{\sfrac{r}{n}}.
    \end{align*}
\end{proof}

The next theorem, originally stated without assuming symmetry, is adapted for
the case of a symmetric initial matrix.
\begin{theorem}[Theorem 5.1 in~\citep{DamleSun19}]
    \label{thm:pert-bound}
    Suppose $\tilde{A} = A + E$ with $A$ symmetric, having eigenvalue
    decomposition
    \(
        A = V \Lambda V^\T + V_{\perp} \Lambda_{\perp} V_{\perp}^\T,
    \)
    where $V \in \Rbb^{n \times r}, V_{\perp} \in \Rbb^{n \times (n - r)}$ have
    orthonormal columns. Moreover, let $\mathsf{gap} := \min\set{\lambda_r -
        \lambda_{r+1}, \mathsf{sep}_{(2, \infty), V_{\perp}}(\Lambda, V_{\perp}
        \Lambda_{\perp} {V_{\perp}}^\T)}$. If $\norm{E}_2 \leq \frac{\mathsf{gap}}{5}$, then the
    leading  invariant subspace of $\tilde{A}$, $\tilde{V}$, satisfies
    \begin{align}
        \begin{aligned}
            \inf_{O \in \Obb_{r}} \norm{\tilde{V} - V O}_{2 \to \infty}
            & \leq 8 \norm{V}_{2 \to \infty}
            \left(\frac{\norm{E}_2}{\lambda_r - \lambda_{r+1}}\right)^2
            + \frac{2 \norm{V_{\perp} {V_{\perp}}^\T E V}_{2 \to
            \infty}}{\mathsf{gap}} \\
            & +
            \frac{4 \norm{V_{\perp} {V_{\perp}}^\T E}_{2 \to \infty}
            \norm{E}_2}{\mathsf{gap} \cdot (\lambda_r - \lambda_{r+1})}.
        \end{aligned}
        \label{eq:pert}
    \end{align}
\end{theorem}

\begin{lemma}[\cite{CapeTangPriebe19}]
    \label{lemma:2inf-subordinate}
    We have
    \begin{align}
        \norm{AB}_{2 \to \infty} &\leq \norm{A}_{2 \to \infty} \norm{B}_2
        \label{eq:submul-right} \\
        \norm{AB}_{2 \to \infty} &\leq \norm{A}_{\infty} \norm{B}_{2 \to \infty}
        \label{eq:submul-left}
    \end{align}
    Moreover, for any matrix $V$ with orthonormal columns, it holds that
    \begin{equation}
        \norm{A V^\T}_{2 \to \infty} = \norm{A}_{2 \to \infty}.
        \label{eq:right-invariance}
    \end{equation}
\end{lemma}

We also prove the following claim, which is used throughout the proof
of~\Cref{prop:2inf-convergence-general} in the next section.
\begin{lemma}
    \label{lemma:2-dist-equiv}
    We have $\inf_{Z \in \Obb_r} \norm{\tilde{V} - V Z}_2 \leq
    \sqrt{2} \dist_2(V, \tilde{V})$.
\end{lemma}
\begin{proof}
    Recall the solution of the orthogonal Procrustes problem, given by
    the SVD of $V^\T \tilde{V}$, $U \Sigma W^\T$. Since $UW^\T \in \Obb_r$,
    with $U^\T U = UU^\T = W^\T W = WW^\T = I_r$, we have
    \begin{align}
        \inf_{Z \in \Obb_r} \norm{\tilde{V} - V Z}_2 &\leq
        \norm{\tilde{V} - V U W^\T}_2 =
        \sqrt{\sup_x \ip{x, (\tilde{V} - VUW^\T)^\T(\tilde{V} - VUW^\T)x}}
        \\
        &= \sqrt{\sup_{x}\ip{x, (I - \tilde{V}^\T V U W^\T - W U^\T
        V^\T \tilde{V} + I) x}} \\
        &\overset{(\sharp)}{=} \sqrt{\sup_{x} \ip{x, 2(I - W \Sigma W^\T)x}} =
        \sqrt{2 \norm{I - W \Sigma W^\T}_2} \\
        &= \sqrt{2} \sqrt{\norm{I - \Sigma}_2} = \sqrt{2} \sqrt{1 -
        \sigma_r(V^\T \tilde{V})} \\
        &\overset{(\natural)}{\leq}
        \sqrt{2} \sqrt{1 - \sigma_r^2(V^\T \tilde{V})}
        = \sqrt{2} \norm{V^\T \tilde{V}}_2,
    \end{align}
    where $(\sharp)$ follows after replacing $V^\T \tilde{V} = U \Sigma
    W^\T$ in the expression and gathering terms, while $(\natural)$ simply
    uses the fact that $\sigma_r(V^\T \tilde{V}) \leq 1$ to upper bound the
    expression inside the square root. Finally, we use the fact
    that:
    \[
        1 - \sigma^2_{\min}(V^\T \tilde{V}) = \norm{V_{\perp}^\T
        \tilde{V}}_2^2 = \dist_2^2(V, \tilde{V}).
    \]
\end{proof}

\section{Omitted proofs}
\subsection{Proof of~\Cref{prop:2inf-convergence-general}}
Starting with the definition of the $2\to\infty$ distance, we have
\begin{align}
    \dist_{2 \to \infty}(Q_t, V) &= \inf_{Z \in \Obb_r}
    \norm{Q_t - V Z}_{2 \to \infty}
    = \inf_{Z \in \Obb_r} \norm{(VV^\T + V_{\perp} V_{\perp}^\T) (Q_t
    - V Z)}_{2 \to \infty} \\
                                 & \overset{(\sharp)}{\leq}
    \sqrt{2} \norm{VV^\T}_{2 \to \infty} \dist_2(Q_t, V) +
    \norm{V_{\perp} V_{\perp}^\T (Q_t - V Z)}_{2 \to \infty}
    \label{eq:step-1-general}
\end{align}
where ($\sharp$) follows from~\Cref{lemma:2inf-subordinate} and the fact
that $\inf_{Z \in \Obb_r} \norm{Q_t - V Z}_2 \leq \sqrt{2} \dist_2(Q_t, V)$. At this
point, note that standard convergence results~\citep{Saad11,GVL13} state
that
\[
    \dist_2(Q_t, V) \leq \left(\frac{\lambda_{r+1}}{\lambda_r}\right)^t
    \frac{d_0}{\sqrt{1 - d_0^2}},
\]
and additionally $\norm{VV^{\T}}_{2\to\infty} \leq \mu \sqrt{\frac{r}{n}}$,
where $\mu$ is the coherence of $V$.

For the remainder, let us first recall a fact from the analysis of subspace
iteration; the $t^{\text{th}}$ iterate $Q_t$ satisfies
\begin{align}
    Q_t R_t &= A^t V^{(0)}, \; \text{ with $R_t$ invertible} \quad
    \Rightarrow
    \quad V_{\perp}^\T Q_t = V_{\perp}^\T A^t V^{(0)} R_t^{-1} =
    \Lambda_{\perp}^t V_{\perp}^\T V^{(0)} R_t^{-1}.
    \label{eq:QtRt}
\end{align}
Then, notice that $V_{\perp}^\T V = 0$ and therefore we can
rewrite the second term in~\eqref{eq:step-1-general} as
\begin{align}
    \norm{V_{\perp} V_{\perp}^\T Q_t}_{2 \to \infty}
    & \overset{(*)}{=}
    \norm{V_{\perp} \Lambda_{\perp}^t V_{\perp}^\T Q_0 R_t^{-1}}_{2 \to
    \infty} \overset{(\flat)}{=} \inf_{Z \in \Obb_r}
    \norm{V_{\perp} \Lambda_{\perp}^t V_{\perp}^\T (Q_0 - V Z)
    R_t^{-1}}_{2 \to\infty} \\
    & \overset{(\natural)}{\leq} \inf_{Z \in \Obb_r}
    C \norm{V_{\perp} V_{\perp}^\T}_{\infty} \lambda_{r + 1}^t
    \norm{(Q_0 - V Z) R_t^{-1}}_{2 \to \infty} \\
    & \leq
    C \norm{V_{\perp} V_{\perp}^\T}_{\infty} \lambda_{r + 1}^t
    \underbrace{\inf_{Z \in \Obb_r} \norm{Q_0 - V Z}_{2 \to \infty}}_{
    = \dist_{2 \to \infty}(Q_0, V)} \norm{R_t^{-1}}_2
    \label{eq:step-2}
\end{align}
where $(*)$ follows from~\cref{eq:QtRt}, $(\flat)$ holds since we can
reintroduce $VZ$ for any $Z$, as $V_{\perp}^\T V = 0$, $(\natural)$ holds
after combining~\cref{eq:submul-left} and
Assumption 1 from the main text, and the last inequality
is~\cref{eq:submul-right}. Notice that
\(
    \norm{R_t^{-1}}_2 = \frac{1}{\sqrt{1 - d_0^2}} \lambda_r^{-t},
\)
by tracing the proof of~\citep[Theorem 8.2.2]{GVL13}. Finally,
by~\Cref{lemma:incoherence-compl}, $\norm{V_{\perp}V_{\perp}^\T}_{\infty}
\leq 1 + \mu \sqrt{r}$.
\qed

\subsection{Proof of~\Cref{prop:2infnorm-convergence-single}}
    For simplicity, let us define $\tilde{V} := \bmx{V & v_{r+1}} \in \Rbb^{n
    \times (r+1)}$ and $\tilde{V}_{\perp}$ for the remaining $n - r - 1$
    eigenvectors forming a basis of $\Rbb^n$. Similarly, let
    $\tilde{\Lambda}_{\perp} = \diag(\lambda_{r+2}, \dots, \lambda_n)$.
    Starting with the definition of the $2\to\infty$ distance, we have
        \begin{align}
            \begin{aligned}
            \dist_{2 \to \infty}(Q_t, V) &= \inf_{Z \in \Obb_r}
            \norm{Q_t - V Z}_{2 \to \infty}
            = \inf_{Z \in \Obb_r} \norm{(VV^\T + V_{\perp} V_{\perp}^\T)
            (Q_t
            - V Z)}_{2 \to \infty} \\
            & \overset{(\sharp)}{\leq}
            \sqrt{2} \norm{VV^\T}_{2 \to \infty} \dist_2(Q_t, V) +
            \norm{V_{\perp} V_{\perp}^\T (Q_t - V Z)}_{2 \to \infty}
            \end{aligned},
            \label{eq:step-1}
        \end{align}
    where ($\sharp$) follows from~\Cref{lemma:2inf-subordinate} in the main
    text and the fact that $\inf_{Z \in \Obb_r} \norm{Q_t - V Z}_2 \leq
    \sqrt{2} \dist_2(Q_t, V)$. Now we may rewrite the second term as
    \begin{align}
        \begin{aligned}
        \norm{(v_{r+1}v_{r+1}^\T +
        \tilde{V}_{\perp}\tilde{V}_{\perp}^\T)
        Q_t}_{2 \to \infty}
        & \leq
        \norm{v_{r+1} v_{r+1}^\T Q_t}_{2 \to \infty} +
        \norm{\tilde{V}_{\perp}\tilde{V}_{\perp}^\T Q_t}_{2 \to
        \infty}\\
        &=
        \norm{v_{r+1} \lambda_{r+1}^t v_{r+1}^\T Q_0 R_t^{-1}}_{2 \to \infty}+
        \norm{\tilde{V}_{\perp} \tilde{\Lambda}^t_{\perp}
\tilde{V}_{\perp}^\T
        Q_0 R_t^{-1}}_{2 \to \infty}.
        \end{aligned}
        \label{eq:step-2-single}
    \end{align}
    Pulling $\lambda_{r + 1}^t$ out of the first norm
    in~\eqref{eq:step-2-single} yields
    \begin{align*}
        \norm{v_{r+1}v_{r+1}^\T (Q_0 - V Z_{\star})}_{2 \to \infty}
        \norm{R_t^{-1}}_2
        \leq \norm{v_{r+1}v_{r+1}^\T}_{\infty} \dist_{2 \to \infty}(Q_0, V)
        \cdot \frac{\lambda_r^{-t}}{\sqrt{1 - d_0^2}},
    \end{align*}
    after using~\Cref{lemma:2inf-subordinate} and the fact that
    $\norm{R_t^{-1}}_2 \leq \frac{\lambda_r^{-t}}{\sqrt{1 - d_0^2}}$, while the
    second norm in~\eqref{eq:step-2-single} can be upper bounded by
    \begin{align*}
        \norm{\cancel{\tilde{V}_{\perp}} \tilde{\Lambda}^t_{\perp}}_2
        \norm{\tilde{V}_{\perp}^\T Q_-}_2 \norm{R_t^{-1}}_2
        = \left( \frac{\lambda_{r+2}}{\lambda_r} \right)^t
        \frac{\dist_{2}(Q_0, \tilde{V})}{\sqrt{1 - d_0^2}},
    \end{align*}
    but as the respective subspaces satisfy $\cV \subset \tilde{\cV}$ we
    have $\dist_2(Q_0, \tilde{V}) \leq \dist_2(Q_0, V)$. Combining
    all the ingredients above completes the proof.
\qed

\subsection{Proof of~\Cref{prop:stopping-criterion}}
The condition on $\norm{E}_2$ combined with the assumption that $Q$ is the
leading invariant subspace of the perturbed matrix $A - EQ^\T$ allows us to
apply~\Cref{thm:pert-bound} for the perturbation  $EQ^\T$, from which we
deduce that the approximate eigenvector matrix $V$ satisfies
\begin{align*}
	\dist_{2 \to \infty}(Q, V) \leq
	8 \norm{V}_{2 \to \infty} \left(
		\frac{\norm{E}_2}{\lambda_r - \lambda_{r+1}}
	\right)^2
    + 2 \frac{\norm{V_{\perp} V_{\perp}^\T E Q^\T V}_{2
     \to \infty}}{\gap}
     + 4 \frac{\norm{V_{\perp} V_{\perp}^\T E}_{2 \to \infty}
         \norm{E}_2}{\gap \cdot (\lambda_r -
     \lambda_{r+1})}
\end{align*}
with the appropriate definition of $\gap$. Using~\Cref{lemma:2inf-subordinate}, we
can upper bound the terms above as
\begin{equation}
    \norm{V_{\perp} V_{\perp}^\T E Q^\T V}_{2 \to \infty} \leq
    \norm{V_{\perp} V_{\perp}^\T}_{\infty} \norm{E Q^\T V}_{2 \to \infty}
    \leq \norm{V_{\perp} V_{\perp}^\T}_{\infty} \norm{E}_{2 \to \infty}
	\underbrace{\norm{Q^\T V}_2}_{\leq 1},
\end{equation}
and similarly for the term $\norm{V_{\perp} V_{\perp}^\T E}_{2 \to \infty}$. \qed

\section{Miscellanea}
\subsection{Eigenvalue localization issues}
\label{appendix:eigenvalue-localization}
We briefly address the issue of when we can safely assume that the approximate
invariant subspace $Q$, utilized in~\Cref{prop:stopping-criterion}, is the
\textbf{leading} invariant subspace of the perturbed matrix $A - EQ^\T$. While
the matrix of Ritz values, $S$, is within $\sqrt{2} \norm{E}_2$ distance of a
set of $r$ eigenvalues of $A$, we do not know whether or not these eigenvalues
correspond to the largest (in magnitude) eigenvalues of $A - EQ^\T$.

In this case, one has to appeal to algorithm-specific arguments. Recall that $A$
has spectral decomposition
\(
    A = V \Lambda V^\T + V_{\perp} \Lambda_{\perp} V_{\perp}^\T,
\)
where $\Lambda$ contains the dominant $r$ eigenvalues. Let $Q_{\perp} \in
\Obb_{n, n - r}$ be orthogonal to the approximate eigenvector matrix $Q \in
\Obb_{n, r}$. Then the following
\begin{align*}
    \bmx{Q^\T \\ Q_{\perp}^\T} (A - EQ^\T)
    \bmx{Q & Q_{\perp}} &=
    \bmx{S & Q^\T (A - EQ^\T) Q_{\perp} \\
    Q_{\perp}^\T Q S &
    Q_{\perp}^\T (A - EQ^\T) Q_{\perp}}
    = \bmx{S & Q^\T A Q_{\perp} \\
    \bm{0} & Q_{\perp}^\T A Q_{\perp}}
\end{align*}
is a Schur decomposition of $A - EQ^\T$, with its eigenvalues being the union
$S \cup \Lambda(Q_{\perp}^\T A Q_{\perp})$ -- the objective becomes showing
that $\norm{\Lambda(Q_{\perp}^\T A Q_{\perp})}_2$ is sufficiently
small, after enough progress of the algorithm. By the variational
characterization of singular values for symmetric matrices, we have
\begin{align}
    \norm{Q_{\perp}^\T A Q_{\perp}}_2 &=
    \sup_{x \in \Sbb^{n - 1}} \abs{
    \ip{x, Q_{\perp}^\T A Q_{\perp} x}} \\
    &=
    \sup_{x \in \Sbb^{n - 1}} \abs{
        \ip{x, Q_{\perp}^\T V \Lambda V^\T Q_{\perp} x}
    +   \ip{x, Q_{\perp}^\T V_{\perp} \Lambda_{\perp} V_{\perp}^\T Q_{\perp} x}} \\
    &\overset{(*)}{\leq}
    \abs{\lambda_1(A)} \norm{Q_{\perp}^\T V}_2^2 +
    \abs{\lambda_{r+1}(A)} \cancelto{\leq 1}{\norm{Q_{\perp}^\T V_{\perp}}}
\end{align}
Therefore, as soon as $\dist_2(V, Q) \leq \sqrt{\varepsilon}$, we know that
$\Lambda(Q_{\perp}^\T A Q_{\perp}) \leq \abs{\lambda_1} \varepsilon +
\abs{\lambda_{r+1}}$; thus when both $\norm{E}_2$ and $\varepsilon$ are small enough,
we can ``match'' $S$ with the leading invariant subspace of $A - EQ^\T$, via
the leading eigenvalues of $A$ itself.

\subsection{Convergence of Procrustes solution}
\label{appendix:procrustes-solution}
Let $V_1, \hat{V}_1$ be a pair of matrices with orthogonal columns. Recall that
the Procrustes solution is the solution to the following matrix nearness problem:
\begin{equation}
    Z_{F} := \argmin_{Z \in \Obb_r} \norm{\hat{V}_1 Z - V_1}_F,
    \label{eq:Z-procrustes}
\end{equation}
for which the solution is available via the SVD of $\hat{V}_1^\T V_1$~\cite{Higham88}.
For the iterates $\set{Q_t}_{t \in \Nbb}$ produced by~\cref{alg:subspace-iteration} in
the main text, notice that
\begin{equation}
    \inf_{Z \in \Obb_r} \norm{Q_t - V Z}_{2 \to \infty} \leq
    \norm{Q_t - V Z_{F}}_{2 \to \infty}
    \leq
    \mu\sqrt{\frac{r}{n}} \norm{Q_t - V Z_F}_{2} +
    \norm{V_{\perp}V_{\perp}^\T Q_t}_{2 \to \infty}.
    \label{eq:procrustes-iterate-decomp}
\end{equation}
For the first term, using the definition of $Z_F$ and choosing
$Z_2 := \argmin_{Z \in \Obb_{r}} \norm{Q_t - VZ}_2$, we may obtain
\begin{align}
    \norm{Q_t - V Z_F}_2 &\leq \norm{Q_t - V Z_F}_F
    \leq \norm{Q_t - V Z_2}_F \\
                         & \overset{(\sharp)}{\leq}
    \sqrt{2r} \cdot \norm{Q_t - V Z_2}_2
    \overset{(\flat)}{\leq} 2 \sqrt{r} \cdot \dist_2(Q_t, V),
\end{align}
where $(\sharp)$ follows by the fact that $\rank(Q_t - V Z_2) \leq 2r$ combined
with norm equivalence, and $(\flat)$ follows from~\cref{lemma:2-dist-equiv}.
Together with the second term in~\cref{eq:procrustes-iterate-decomp}, these can
be analyzed as in the proofs of~\Cref{prop:2inf-convergence-general,prop:2infnorm-convergence-single}.
.

\subsection{Convergence \textit{without} Assumption~\ref{asm:norm-bound}}
\label{appendix:preliminary-results}
Here, we provide a proof showing that the convergence of subspace iteration
w.r.t. the $2 \to \infty$ norm improves upon the spectral norm results without
the need for Assumption 1 from the main text on the data matrix's eigenspaces. We show this by
studying a ``worse-case'' version of $A$, $\tilde{A}$, instead; given $A =
V \Lambda V^\T + V_{\perp} \Lambda_{\perp} V_{\perp}^\T$, we define $\tilde{A}$ as
\begin{equation}
    \tilde{A} := V \Lambda V^\T + \lambda_{r+2}(A) \cdot V_{\perp} {V_{\perp}}^\T.
\end{equation}
In the forthcoming proof, we denote $\tilde{\Lambda}_{\perp} := \lambda_{r+2}(A) \cdot I_{n-r}$.
The Proposition below gives an improved rate compared to the analysis w.r.t.
spectral norm convergence, albeit for a limited set of spectra.

\begin{myprop}
	\label{prop:subspace-iter-2inf}
	The iterates $\set{Q_t}_{t \in [T]}$ produced by Algorithm~\ref{alg:subspace-iteration}
    in the main text with initial guess $V^{(0)}$ satisfy
	\begin{align}
	\label{eq:subspace-iter-2-inf-dist}
    \begin{aligned}
    \dist_{2 \to \infty}(Q_t, V)
    & \leq 3 \frac{1 + \mu
    \sqrt{r}}{\sqrt{1 - d_0^2}}\left(
	\frac{\lambda_{r+2}}{\lambda_r} \right)^t \cdot
	\dist_{2 \to \infty}(V^{(0)}, V) \\
	& \quad +
	\mu \sqrt{\frac{r}{n}}
	\left(\frac{\lambda_{r + 1}}{\lambda_r}\right)^t
	\cdot \tan(\theta_0)
    + \max\set{ \frac{\lambda_{r + 1}^t - \lambda_{r +
	2}^t}{\lambda_r^t},
		\frac{\lambda_{r + 2}^t - \lambda_n^t}{\lambda_r^t}}
	\cdot \tan(\theta_0),
    \end{aligned}
	\end{align}
    where $\tan(\theta_0) := \frac{d_0}{\sqrt{1 - d_0^2}}, \; d_0 := \dist_2(Q_0, V)$.
\end{myprop}
\begin{proof}
Let us introduce some notation to be used in the proof; given the true subspace
$Q$, we write $\dist_{\norm{\cdot}, \perp}(A, B) :=
\dist_{\norm{\cdot}}(V_{\perp} V_{\perp}^\T A, V_{\perp} V_{\perp}^\T
B)$. By splitting up $Q_t$ into its projections to $V$ and $V_{\perp}$
respectively, we can upper bound the desired distance in the following way:
\begin{align}
	\dist_{2, \infty}(Q_t, V)
    &= \inf_{Z} \norm{Q_t - V Z}_{2 \to \infty} =
	\inf_{Z} \norm{(VV^\T + V_{\perp} V_{\perp}^\T) Q_t - VZ}_{2 \to
	\infty} \notag \\
    & \leq
	\inf_{Z} \norm{VV^\T(Q_t - V Z)}_{2 \to \infty}
	+ \norm{V_{\perp} V_{\perp}^\T Q_t}_{2 \to \infty} \\
    & \leq \norm{V}_{2 \to \infty} \dist_2(Q_t, V) + \dist_{2 \to \infty, \perp}(Q_t, V)
	\label{eq:init-decomp}
\end{align}
since $V_{\perp}^\T V = 0$. The first term in~\eqref{eq:init-decomp} is upper bounded by
\[
	\mu \sqrt{\frac{r}{n}} \left(\frac{\lambda_{r+1}}{\lambda_r}\right)^t
	\tan (\theta_0),
\]
(where $\mu$ is the coherence of $V$),
which is known from the standard convergence analysis of
Algorithm~\ref{alg:subspace-iteration} measured in the spectral norm.
In addition, using the triangle inequality for the second term
in~\eqref{eq:init-decomp}, we can further upper bound
\begin{align}
    \begin{aligned}
    \norm{V_{\perp} V_{\perp}^\T Q_t}_{2 \to \infty} &=
	\dist_{2\to\infty,\perp}(Q_t, V)
	\leq \dist_{2\to\infty,\perp}(Q_t, \tilde{Q}_t) +
	\dist_{2\to\infty,\perp}(\tilde{Q}_t, V)
    \end{aligned}
	\label{eq:dist-decomp}
\end{align}
where $\tilde{Q}_t$ is the aforementioned ``ghost'' iterate resulting from
applying Algorithm~\ref{alg:subspace-iteration} to the matrix $\tilde{A}$,
which is defined as
\begin{equation}
	\tilde{A} := \bmx{V & V_{\perp}} \bmx{\Lambda & 0 \\
    0 & \lambda_{r+2}(A) I_{n - r}} \bmx{V & V_{\perp}}^\T
	\label{eq:pert-A},
\end{equation}
and obviously $Q_0 = \tilde{Q}_0 := V^{(0)}$.
In the forthcoming steps, we bound each distance above separately. For
the second term in~\eqref{eq:dist-decomp}, we have:
\begin{lemma}
	\label{lemma:subspace-iter-qtilde}
    The iterates $\set{\tilde{Q}_t}_{t \in [T]}$ produced by
    Algorithm~\ref{alg:subspace-iteration} when applied to
    $\tilde{A}$, as defined in~\eqref{eq:pert-A}, satisfy
    \begin{flalign}
        \begin{aligned}
		\dist_{2 \to \infty, \perp}(\tilde{Q}_t, V) &\leq
        \left( \frac{\lambda_{r+2}}{\lambda_r} \right)^t
        \frac{1}{\sqrt{1 - d_0^2}}
            (1 + \mu \sqrt{r}) \cdot \dist_{2 \to \infty}(Q_0, V) &
        \end{aligned}
    \end{flalign}
    where $d_0 := \dist_{2}(Q_0, V)$, and $\mu$ is the coherence of $V$.
\end{lemma}
\begin{proof}
We build heavily on the proof of the analogous convergence result for the
spectral norm given in~\citep{GVL13}. First, observe that $\tilde{A}$ has the
same eigenvectors as $A$ and same first $r$ as well as last $n - r - 1$
eigenvalues.

From the proof of~\citep[Theorem 8.2.2]{GVL13}, we know that $\tilde{Q}_t
\tilde{R}_t = \tilde{A}^t V^{(0)}$,
with $\tilde{R}_t$ invertible and satisfying
\begin{equation}
    \norm{\tilde{R}_t^{-1}}_2 \leq \frac{\lambda_r^{-t}}{\sqrt{1 - d_0^2}},
    \quad d_0 := \dist_{2}(V^{(0)}, V),
    \label{eq:Rinv-norm}
\end{equation}
Then we have
\begin{align}
    \begin{aligned}
        \dist_{2 \to \infty,\perp}(\tilde{Q}_t, V) &=
	\norm{V_{\perp} V_{\perp}^\T \tilde{Q}_t}_{2 \to \infty} =
	\norm{V_{\perp} \tilde{\Lambda}_{\perp}^\T V_{\perp}^\T V^{(0)}
	\tilde{R}_t^{-1}}_{2 \to \infty} \\
                                                   &\leq
	\inf_{Z} \norm{V_{\perp} \tilde{\Lambda}_{\perp}^\T V_{\perp}^\T
	(V^{(0)} - V Z)}_{2 \to \infty} \norm{\tilde{R}_t^{-1}}_{2} \\
    & \overset{(\sharp)}{\leq}
	\left( \frac{\lambda_{r+2}}{\lambda_r}\right)^t
	\norm{V_{\perp} V_{\perp}^\T}_{\infty}
	\dist_{2 \to \infty}(V^{(0)}, V) \frac{1}{\sqrt{1 - d_0^2}}
    \end{aligned}
\end{align}
where the second step of the proof uses~\eqref{eq:submul-right} and
$(\sharp)$ uses~\Cref{eq:submul-left,eq:Rinv-norm}.
Finally, an appeal to~\Cref{lemma:incoherence-compl} yields the desired
expression.
\end{proof}

For the first term in~\eqref{eq:dist-decomp}, we follow a similar approach and
write (for $Z$ attaining the infimum in the definition of the subspace
distance):
\begin{align*}
	\dist_{2 \to \infty, \perp}(Q_t, \tilde{Q}_t) =
	\norm{V_{\perp} V_{\perp}^\T (Q_t - \tilde{Q}_t Z)}_{2 \to
	\infty}
												  =
	\norm{V_{\perp} \Lambda_{\perp}^t V_{\perp}^\T R_t^{-1} -
		  V_{\perp} \tilde{\Lambda}_{\perp}^t V_{\perp}^\T \tilde{R}_t^{-1}Z}_{2 \to \infty},
\end{align*}
where again we recall from the proof of~\citep[Theorem 8.2.2]{GVL13} that
\begin{align*}
	V_{\perp}^\T Q_t = \Lambda_{\perp}^t V_{\perp}^\T V^{(0)}
	R^{-1}_{t}, \; V_{\perp}^\T \tilde{Q}_t = \lambda_{r + 2}^t V_{\perp}^\T
	V^{(0)} \tilde{R}_{t}^{-1},
\end{align*}
as in the proof of~\Cref{lemma:subspace-iter-qtilde}.
Consequently, we can write $\tilde{R}_t^{-1} Z =
R_t^{-1} + (\tilde{R}_t^{-1} Z - R_t^{-1})$ and substitute above to obtain
\begin{align}
	\begin{aligned}
	\norm{V_{\perp} V_{\perp}^\T (Q_t - \tilde{Q}_t Z)}_{2 \to \infty}
    & \leq
    \norm{V_{\perp} \left(\Lambda_{\perp}^t - \lambda_{r + 2}^t I_{n - r}
		\right) V_{\perp}^\T V^{(0)} R_t^{-1}}_{2 \to \infty} \\
	& \quad + \lambda_{r + 2}^t \norm{V_{\perp} V_{\perp}^\T V^{(0)}
		(\tilde{R}_t^{-1} Z - R_{t}^{-1})}_{2 \to \infty},
	\end{aligned}
	\label{eq:Qtild-decomp}
\end{align}
after appropriate rearrangements and the triangle inequality.
Rewriting $V_{\perp}^\T V^{(0)} = V_{\perp}^\T  (V^{(0)} - V Z_{2 \to
\infty})$ yields
\begin{align}
	\begin{aligned}
	\norm{V_{\perp} V_{\perp}^\T V^{(0)} (\tilde{R}_t^{-1} - R_{t}^{-1})}_{2
	\to \infty} \leq
	\norm{V_{\perp} V_{\perp}^\T}_{\infty} \dist_{2 \to \infty}(V^{(0)}, V)
	\left( \norm{\tilde{R}_t^{-1}}_2 + \norm{R_t^{-1}}_2 \right) \\
	\leq 2 \left(1 + \mu \sqrt{r}\right)
	\frac{1}{\lambda_r^t} \frac{\dist_{2,\infty}(V^{(0)}, V)}{\sqrt{1 - d_0^2}},
	\end{aligned}
	\label{eq:Qtild-term-2}
\end{align}
since $R_t^{-1}$ and $\tilde{R}_t^{-1}$ both satisfy~\eqref{eq:Rinv-norm}
as $A$ and $\tilde{A}$ have the same first $r$ eigenvalues. Finally
\begin{align}
	\norm{\Lambda_{\perp}^t - \lambda_{r + 2}^t I_{n-r}}_2 &=
	\max\set{\abs{\lambda_{r + 1}^t - \lambda_{r+2}^t},
		\abs{\lambda_{r+2}^t - \lambda_n^t}},
    \label{eq:lambda-term}
\end{align}
which we may use to bound the first term in~\eqref{eq:Qtild-decomp} by noticing
\begin{align}
    \begin{aligned}
    \norm{V_{\perp}(\Lambda_{\perp}^t - \lambda_{r + 2}^t I) V_{\perp}^\T
    V^{(0)} R_t^{-1}}_{2 \to \infty}
    \leq
    \norm{V_{\perp}}_{2 \to \infty}
    \norm{\Lambda_{\perp}^t - \lambda_{r + 2}^t}_2
    \norm{V_{\perp}^\T V^{(0)}}_2 \norm{R_t^{-1}}_2
    \end{aligned}
	\label{eq:Qtild-term-1}
\end{align}
The proof follows by
combining~\Cref{eq:Qtild-term-1,eq:Qtild-term-2,eq:lambda-term}, the fact that
$d_0 = \norm{V_{\perp}^\T V^{(0)}}_2$, and
\Cref{lemma:subspace-iter-qtilde}.
\end{proof}

\section{Reproducibility}
We provide an open-source implementation of the algorithms and all experiments
in \texttt{Julia} in the following repository: \url{https://github.com/VHarisop/entrywise-convergence}.
The experiments were run in a machine running Manjaro Linux
with $16$ GB of RAM and Intel\textregistered Core\texttrademark~i7-7700 CPU @ 3.60 GHz, using Julia version 1.1.

\end{document}